\input ./style/arxiv-general.cfg
\documentclass[MSNbibl,number,citesort,seceqn,secthm,dvips]{arxbj}
\makeatletter
   \@ifpackageloaded{graphicx}{}{\usepackage{graphicx}}
\makeatother

\aid{0}
\volume{22}
\issue{1}
\pubyear{2016}
\firstpage{589}
\lastpage{614}
\doi{10.3150/14-BEJ669} 
\docsubty{FLA}

\makeatletter

\newcommand{\rrvert}{\vert}
\newcommand{\llvert}{\vert}
\newcommand{\cal}{\mathcal}
\newremark{rem}{Remark}[section]
\makeatother

\begin{document}
\begin{frontmatter}

\title{On differentiability of implicitly defined function in semi-parametric profile likelihood~estimation}
\runtitle{On differentiability of implicitly defined function}

\begin{aug}
\author{\inits{Y.}\fnms{Yuichi}~\snm{Hirose}\corref{}\ead
[label=e1]{Yuichi.Hirose@msor.vuw.ac.nz}}
\address{School of Mathematics, Statistics and Operations Research,
Victoria University of Wellington,
New~Zealand. \printead{e1}}
\end{aug}

\received{\smonth{8} \syear{2013}}
\revised{\smonth{4} \syear{2014}}

%
\begin{abstract}
In this paper, we study the differentiability of implicitly defined functions
which we encounter in the profile likelihood estimation of parameters
in semi-parametric models.
Scott and Wild (\textit{Biometrika} \textbf{84} (1997)  57--71;
\textit{J. Statist. Plann. Inference} \textbf{96} (2001) 3--27)
and Murphy and van der Vaart
(\textit{J. Amer. Statist. Assoc.} \textbf{95} (2000) 449--485) developed methodologies
that can avoid dealing with such implicitly defined functions by
parametrizing parameters in the profile likelihood
and using an approximate least favorable submodel in semi-parametric models.
Our result shows applicability of an alternative approach presented in
Hirose (\textit{Ann. Inst. Statist. Math.} \textbf{63} (2011)
1247--1275)
which uses the direct expansion of the profile likelihood.
\end{abstract}

%
\end{frontmatter}

\section{Introduction}\label{sec1}

Consider a general semi-parametric model
\[
{\cal P}=\bigl\{p_{\theta,\eta}(x)\dvt \theta\in\Theta,\eta\in H \bigr\},
\]
where $p_{\theta,\eta}(x)$ is a density function on the sample space
${\cal X}$ which depends on a finite-dimensional parameter $\theta$
and an infinite-dimensional parameter $\eta$.
We assume that the set $\Theta$ of the parameter $\theta$ is an open
subset of $R^d$ and the set $H$ is a convex subset of a Banach space
${\cal B}$.

Once observations $X_1,\ldots,X_n$ are generated from the model, the
log-likelihood is given by
%
\begin{equation}
\ell_n(\theta,\eta)=n^{-1}\sum
_{i=1}^n\log p_{\theta,\eta}(X_i) =
\int \log p_{\theta,\eta}(x) \,\mathrm{d}F_n(x),
\end{equation}
where $F_n$ is the empirical c.d.f. based on the observations.
In the profile likelihood approach, we find a function $\eta_{\theta
,F}$ of the parameter $\theta$ and a c.d.f. $F$
as the maximizer of the log-likelihood given~$\theta$ such that
%
\begin{equation}
\eta_{\theta,F_n}=\mathop{\arg\max}_{\eta}\int\log p_{\theta,\eta}(x)
\,\mathrm{d}F_n(x).
\end{equation}

Then the profile (log)-likelihood is given by
%
\begin{equation}
\int\log p_{\theta,\eta_{\theta,F_n}}(x) \,\mathrm{d}F_n(x).
\end{equation}

In this paper, we consider the situation when the function $\eta_{\theta
,F}$ is given as the solution to the operator equation of the form
%
\begin{equation}
\label{ImplicitEquation}
\eta= \Psi_{\theta,F}(\eta).
\end{equation}

Murphy, Rossini and van der Vaart \cite{MRV97} encountered this type of implicitly defined function in
their maximum likelihood estimation problem in the proportional odds model.
According to them, ``because $\hat{H}_{\beta}$ is not an explicit
function of $\beta$,
we are unable to differentiate the profile log-likelihood explicitly in
$\beta$ to form an estimator of $\Sigma$''
(here $\hat{H}_{\beta}$ is the maximizer of the log-likelihood $\ell
_n(\beta,H)$ given $\beta$,
$H$ is the baseline odds of failure and $\Sigma$ is the efficient information).
The authors (Murphy, Rossini and van der Vaart \cite{MRV97}) used a numerical approximation to the problem.
In the first example (Example~1) given below, we present a modified
version of the proportional odds model
and give an example of implicitly defined function there.

Scott and Wild \cite{SW97,SW01} also encountered implicitly defined functions in their
estimation problem
with data from various outcome-dependent sampling design. They proposed
a method of re-parametrization of profile-likelihood
so that the log-likelihood is an explicitly defined function in terms
of the parameters in the re-parametrized model.
Their estimators turned out to be efficient and Hirose and Lee \cite{HL12} showed
conditions under which re-parametrization gives efficient estimation
in a context of multiple-sample semi-parametric model.

Another way to avoid dealing with implicitly defined functions is
developed by Murphy and van der Vaart \cite{MV00}.
The paper proved the efficiency of profile likelihood estimation
by introducing an approximate least favorable sub-model
to express the upper
and lower bounds for the profile log-likelihood.
Since these two bounds have the same expression for the asymptotic
expansion, so does the one for the profile log-likelihood.
The advantage of the approach is that it does not need to deal with
implicitly defined functions which we discussed in the current paper.
Disadvantage of Murphy and van der Vaart \cite{MV00} are
(1) it needs to find an approximate least favorable submodel in each
example 
which may be difficult to find in some cases;
(2) no-bias condition (equation (3.4) in Murphy and van der Vaart \cite{MV00}) is assumed in the
main theorem and it needs to be verified in examples to which the main
theorem is applied.
In their ``Discussion'', they commented
``It appears difficult to derive good approximations to a least
favorable path for such models,
and given such approximation it is unclear how one would verify the
no-bias condition''.

Hirose \cite{Hirose10} used direct asymptotic expansion of the profile
likelihood to show the efficiency of the profile likelihood estimator.
The result in the paper (Theorem~1 in Hirose \cite{Hirose10}) does not assume
the no-bias condition and, under the assumptions given there,
the no-bias condition (equation~(4) in Hirose \cite{Hirose10}) is proved
(therefore, verification of the no-bias condition is not required in
examples). 
In the approach, we cannot avoid dealing with implicitly defined
functions of the form given in (\ref{ImplicitEquation}) in some applications.
The purpose of this paper is to study the properties of these function
such as differentiability
so that the method in Hirose \cite{Hirose10} is applicable to those applications.
The results in Hirose \cite{Hirose10} are summarized in Section~\ref{sec6}.

In Section~\ref{sec2}, we give examples of implicitly defined functions.
The main results are presented in Section~\ref{sec3}.
In Sections~\ref{sec4} and~\ref{sec5}, the main results are applied to
the examples.
In Section~\ref{sec61}, we demonstrate  how the result of the paper
(the differentiability of implicitly defined functions in
semi-parametric models)
can be applied in a context of asymptotic linear expansion of the
maximum profile likelihood estimator in a semi-parametric model.

\section{Examples}\label{sec2}

\subsection{Example~1 (semi-parametric proportional odds model)}\label{sec21}

The original asymptotic theory for maximum likelihood estimator in the
semi-parametric proportional odds model is developed in Murphy, Rossini and van der
Vaart \cite{MRV97}.
We present a modified version of the model in Kosorok \cite{Kosorok08}.

In this model, we observe $X=(U,\delta,Z)$, where $U=T \wedge C$,
$\delta=1_{\{U=T\}}$,
$Z \in R^d$ is a covariate vector, $T$ is a failure time and $C$ is a
right censoring time.
We assume $C$ and $T$ are independent given $Z$.

The proportional odds regression model is specified by the survival
function of $T$ given $Z$ of the form
\[
S(t|Z)=\frac{1}{1+\mathrm{e}^{\beta'Z}A(t)},
\]
where $A(t)$ is nondecreasing function on $[0,\tau]$ with $A(0)=0$.
$\tau$ is the limit of censoring distribution such that $P(C>\tau)=0$
and $P(C=\tau)>0$.
The distribution of $Z$ and $C$ are uninformative of $S$ and $\operatorname{var} Z$ is positive definite.

Define the counting process $N(t)=\delta1_{\{U \leq t\}}$ and at risk
process $Y(t)=1_{\{U \geq t\}}$.
We assume $P\{\delta Y(t)=1\}>0$ for each $t \in[0,\tau]$.

Let $F_n$ be the empirical process for i.i.d. observation $(U_i,\delta
_i,Z_i)$, $i=1,\ldots,n$.
Then the log-likelihood on page 292 in Kosorok \cite{Kosorok08} can be written as
\[
\label{loglikeex1}
\ell_n(\beta,A)
 =  \int\bigl\{\delta\bigl(
\beta'Z + \log a(U)\bigr)-(1+\delta)\log\bigl(1+\mathrm{e}^{\beta
'Z}A(U)
\bigr)\bigr\} \,\mathrm{d}F_n,
\]
where $a(t)=\mathrm{d}A(t)/\mathrm{d}t$.

Consider one-dimensional sub-models for $A$ defined by the map
\[
t \rightarrow A_t(u)=\int_0^u
\bigl(1+th(s)\bigr)\,\mathrm{d}A(s),
\]
where $h(s)$ is an arbitrary total variation bounded cadlag function on
$[0,\tau]$.
By differentiating the log-likelihood function $\ell_n(\beta,A_t)$
with respect to $t$ at $t=0$, we obtain the score operator
\[
B_n(\beta,A) (h)  =  \frac{\mathrm{d}}{\mathrm{d}t} \bigg|_{t=0}
\ell_n(\beta,A_t) = \int \biggl\{ \delta h(U)-(1+\delta)
\frac{\mathrm{e}^{\beta'Z}\int_0^{U}h(u)\,\mathrm{d}A(u)}{1+\mathrm{e}^{\beta'Z}A(U)} \biggr\} \,\mathrm{d}F_n.
\]

Choose $h(u)=1_{\{u \leq t\}}$, then
\[
B_n(\beta,A) (h)  =  \int N(t) \,\mathrm{d}F_n - \int \biggl\{
\int_0^{U}W(u;\beta,A)\,\mathrm{d}A(u) \biggr
\}\,\mathrm{d}F_n,
\]
where
$N(t)$ 
and
$Y(t)$
are defined above and
%
\begin{equation}
\label{Wt}
W(u;\beta,A)=\frac{(1+\delta)\mathrm{e}^{\beta'Z}Y(u)}{1+\mathrm{e}^{\beta'Z}A(U)}.
\end{equation}
The solution $\hat{A}_{\beta,F_n}$ to the equation $B_n(\beta,A)(h)=0$
is of the form
%
\begin{equation}
\label{hatAbeta}
\hat{A}_{\beta,F_n}(u)=\int_0^u
\frac{E_{F_n}\, \mathrm{d}N(s) }{E_{F_n} W(s;\beta
,\hat{A}_{\beta,F_n})},
\end{equation}
where $E_{F_n}\,\mathrm{d}N(s)= \int \mathrm{d}N(s) \,\mathrm{d}F_n$ and $E_{F_n} W(s;\beta,\hat
{A}_{\beta,F_n})=\int W(s;\beta,\hat{A}_{\beta,F_n})\,\mathrm{d}F_n$.

Let $F$ be a generic notation for the c.d.f., and if we let
%
\begin{equation}
\label{PsiPodds}
\Psi_{\beta,F}(A)=\int_0^u
\frac{E_F \,\mathrm{d}N(s)}{E_F W(s;\beta,A)},
\end{equation}
then (\ref{hatAbeta}) is a solution to the operator equation $A=\Psi
_{\beta,F_n}(A)$, here
$E_{F}\,\mathrm{d}N(s)= \int \mathrm{d}N(s) \,\mathrm{d}F$ and $E_{F} W(s;\beta,\hat{A}_{\beta
,F})=\int W(s;\beta,\hat{A}_{\beta,F})\,\mathrm{d}F$.
More detailed treatment of this example can be found in \cite{Kosorok08}, Section~15.3, pages 291--303.
We continue this example in Section~\ref{sec4}.

\subsection{Example~2 (continuous outcome with missing data)}\label{sec22}
This example is studied in Weaver and Zhou \cite{WZ05} and Song, Zhou and Kosorok \cite{SZK09}.
Suppose the underlying data generating process on the sample space
${\cal Y}\times{\cal X}$ is a model
%
\begin{equation}
\label{modelQ}
{\cal Q}=\bigl\{p(y,x;\theta)=f(y|x;\theta)g(x)\dvt \theta\in\Theta,
g \in {\cal G} \bigr\}.
\end{equation}
Here, $f(y|x;\theta)$ is a conditional density of $Y$ given $X$ which
depends on a finite-dimensional parameter $\theta$,
$g(x)$ is an unspecified density of $X$ which is an
infinite-dimensional nuisance parameter.
We assume the set $\Theta\subset R^d$ is an open set containing a
neighborhood of the true value $\theta_0$
and ${\cal G}$ is the set of density function of $x$ containing the
true value $g_0(x)$.
We assume the variable $Y$ is a continuous variable.

We consider a situation when there are samples for which we observe
complete observation $(Y,X)$ and for which we observe only $Y$.
Let $R_i$ be the indicator variable for the $i$th observation defined by
\[
R_i=\left\{\begin{array}{l@{\qquad}l} 1, &\mbox{if $X_i$ is observed,}
\\
2, & \mbox{if $X_i$ is not observed.}\end{array}\right.
\]
Then the index set for the complete observations is $V=\{i\dvt R_i=1\}$ and
the index set for the incomplete observations is $\overline{V}=\{
i\dvt R_i=2\}$.
(In the paper  Song, Zhou and Kosorok \cite{SZK09} \mbox{$R_i=0$} was used for subjects $X_i$ is not
observed.)
Let $n_V=|V|$, $n_{\overline{V}}=|\overline{V}|$ be the total number of
complete observations and incomplete observations, respectively.


Weaver and Zhou \cite{WZ05} and Song, Zhou and Kosorok \cite{SZK09} consider the likelihood of the\vspace*{-2pt} form
%
\begin{equation}
\label{LikelihoodExt}
L_n(\theta,g)  =  \prod_{i \in V}
\bigl\{f(Y_i|X_i;\theta)g(X_i) \bigr\}
\prod_{i \in\overline{V}} f_Y(Y_i;
\theta,g),
\end{equation}
where\vspace*{-3pt}
%
\begin{equation}
\label{fY} f_Y(y;\theta,g)= \int_{\cal X} f(y|x;
\theta)g(x)\,\mathrm{d}x.
\end{equation}


The log-likelihood, the $1/n$ times log of (\ref{LikelihoodExt}) is
\[
\ell_n(\theta,g)  =  \frac{n_V}{n} \frac{1}{n_V}\sum
_{i \in V} \bigl\{\log f(y_{i}|x_{i};
\theta)+\log g(x_{i}) \bigr\} +\frac{n_{\overline{V}}}{n} \frac{1}{n_{\overline{V}}}\sum
_{i \in
\overline{V}} \log f_Y(y_{i};
\theta,g).
\]

For the proof in the later part of the paper, we introduce notation:
let $F_{1n}$ and $F_{2n}$ be the empirical c.d.f.s based on the samples in
$V$ and $\overline{V}$, respectively;
denote $w_{1n}=n_{V}/n$, $w_{2n}=n_{\overline{V}}/n$
and let $F_n=\sum_{s=1}^2w_{sn}F_{sn}$ be the empirical c.d.f. for the
combined samples in $V \cup\overline{V}$.

Then the log-likelihood can be expressed as
\[
\ell_n(\theta,g) 
 =  w_{1n} \int \bigl\{\log
f(y|x;\theta)+\log g(x) \bigr\}\,\mathrm{d}F_{1n} +w_{2n} \int\log
f_Y(y;\theta,g)\,\mathrm{d}F_{2n}.
\]

To find the maximizer of $\ell_n(\theta,g)$, we treat $g(x)$ as
probability mass function on the observed values $\{x_i\dvt i \in V\}$.
Denote $g_i=g(x_i)$, $i \in V$.
The derivative of the log-likelihood with respect to~$g_i$ is
\[
\frac{\partial}{\partial g_i } \ell_n(\theta,g)  =  w_{1n}
\frac{\int1_{\{x=x_i\}}\,\mathrm{d}F_{1n}}{g_i} + w_{2n} \int \frac{
f(y|x_i;\theta)}{f_Y(y;\theta,g)} \,\mathrm{d}F_{2n},
\]
here, for the discrete $g$, $f_Y(y;\theta,g)=\sum_{i \in V}
f(y|x_i;\theta)g_i$.

Let $\lambda$ be a Lagrange multiplier to account for $\sum_{i \in V} g_i=1$.
Set $\frac{\partial}{\partial g_i } \ell_n(\theta,g) + \lambda=0$.
Multiply by $g_i$ and sum over $i \in V$ to get $w_{1n}+w_{2n}+\lambda=0$.
Therefore, $\lambda=-(w_{1n}+w_{2n})=-1$\vspace*{1pt} and
$\frac{\partial}{\partial g_i }\ell_n(\theta,g) -1=0$.
By rearranging this equation, we obtain
\[
\hat{g}_i 
 =
\frac{w_{1n} \int1_{\{x=x_i\}}\,\mathrm{d}F_{1n} }{1 - w_{2n} \int {
f(y|x_i;\theta)}/{f_Y(y;\theta,\hat{g})} \,\mathrm{d}F_{2n}}.
\]
This is exactly equation (3) in Song, Zhou and Kosorok \cite{SZK09}.
Since the $\hat{g}_i$ is a function of $\theta$ and $F_n=\sum_{s=1}^2w_{sn}F_{sn}$, it can be written as
%
\begin{equation}
\label{hatg}
\hat{g}_{\theta,F_n}(x_i) =\frac{ w_{1n}  (\partial_x \int \mathrm{d}F_{1n} )(x_i)}{1 - w_{2n}
\int { f(y|x_i;\theta)}/{f_Y(y;\theta,\hat{g}_{\theta,F_n})}
\,\mathrm{d}F_{2n}},\qquad  i \in V,
\end{equation}
where $\partial_x=\frac{\partial}{\partial x}$ (see Note  below
for the notation $\partial_x \int \mathrm{d}F_{1}$).
This is a solution to the equation $g=\Psi_{\theta,F_n}(g)$ with
\[
\Psi_{\theta,F}(g)=\frac{ w_{1} \partial_x \int \mathrm{d}F_{1} }{1 - w_{2}
\int { f(y|x;\theta)}/{f_Y(y;\theta,g)} \,\mathrm{d}F_{2}},
\]
here $F=\sum_{s=1}^2w_{s}F_{s}$.
We continue this example in Sections~\ref{sec5} and \ref{sec61}.

\textit{Note (Comment on the notation $\partial_x \int\mathrm{d}F_{1}$)}.
Let us denote $\partial_x=\frac{\partial}{\partial x}$. The Heaviside
step function $H(x)=1_{\{x \geq0\}}$ and the Dirac delta function
$\delta(x)$ are related by
$\partial_x H(x)=\delta(x)$.
Using this, for the joint empirical c.d.f. $F_{n}(x,y) =\frac{1}{n}\sum_{i=1}^n H(x-x_i) H(y-y_i)$,
we have
\[
\biggl(\partial_x \int \mathrm{d}F_n \biggr) (x)=
\frac{1}{n}\sum_{i=1}^n \delta
(x-x_i) \int \mathrm{d}H(y-y_i)=\frac{1}{n}\sum
_{i=1}^n \delta(x-x_i),
\]
where we used $\int \mathrm{d}H(y-y_i)=1$ (since the integral is over all $y$).
For the continuous case, joint c.d.f. $F(x,y)$ and marginal p.d.f. $f(x)$
are related by
$(\partial_x \int \mathrm{d}F)(x) =f(x)$.
This justifies the notation $\partial_x \int \mathrm{d}F_{1}$ for both
continuous and empirical c.d.f.s.

\section{Main results}\label{sec3}

In this section, we show the differentiability of implicitly defined
function which 
is given as a solution to the operator equation (\ref{ImplicitEquation}).

As we stated in the \hyperref[sec1]{Introduction}, we consider a general
semi-parametric model
\[
{\cal P}=\bigl\{p_{\theta,\eta}(x)\dvt\theta\in\Theta,\eta\in H \bigr\},
\]
where $p_{\theta,\eta}(x)$ is a density function on the sample space
${\cal X}$ which depends on a finite-dimensional parameter $\theta$
and an infinite-dimensional parameter $\eta$.
We assume that the set $\Theta$ of the parameter~$\theta$ is an open
subset of $R^d$ and the set $H$ is a convex set in a Banach space
${\cal B,}$
which we may assume the closed linear span of $H$.


\textit{Definition (Hadamard differentiability)}. Suppose $X$ and $Y$ are
two normed linear spaces and let $T \subset X$.
We say that a map $\psi\dvtx T \rightarrow Y$ is Hadamard differentiable at
$x \in T$
if there is a~continuous linear map $\mathrm{d}\psi(x)\dvtx X\rightarrow Y$ such that
%
\begin{equation}
\label{Hadmard}
t^{-1}\bigl\{\psi(x_t)-\psi(x)\bigr\}
\rightarrow \mathrm{d}\psi(x)h \qquad \mbox{as } t \downarrow 0
\end{equation}
for any map $t \rightarrow x_t$ with $x_{t=0}=x$ and $t^{-1}(x_t-x)
\rightarrow h \in X$ as $t \downarrow0$.
The map $\mathrm{d}\psi(x)$ is called the Hadamard derivative of $\psi$ at $x$,
and is continuous in $x$
(for reference, see Gill \cite{Gill89} and Shapiro~\cite{Shapiro90}).

We denote the second derivative of $\psi$ in the sense of Hadamard by
$\mathrm{d}^2\psi(x)$.
The usual first and second derivative of a parametric function $\psi
(x)$, $x \in R^d$, are denoted by $\dot{\psi}$ and $\ddot{\psi}$.

\textit{Note on Hadamard differentiability}. The above form of definition
of the Hadamard differentiability is due to Fr\'{e}chet in 1937.
M. Sova showed the equivalence of the Hadamard differentiability and
the compact differentiability in metrizable linear spaces
(Averbukh and Smolyanov \cite{AS68}).
Because of the equivalence, some authors use compact differentiability
as definition of Hadamard differentiability
(Gill \cite{Gill89}, van der Vaart and Wellner \cite{VW96}, Bickel, Klaassen, Ritov and
Wellner \cite{BKRW93}).
In this paper, we use the definition of Hadamard differentiability
given by Fr\'{e}chet.

In addition to the Hadamard differentiability of functions, in Theorem~\ref{mainthm} below, we assume the following condition.

\textit{Additional condition.}
We say a Hadamard differentiable map $\psi(x)$ satisfies the additional
condition at $x$, if,
for each path $x_t$ in some neighborhood of $x$,
there is a bounded and linear map $h \rightarrow \mathrm{d}\psi^*_th$ such that
the equality
%
\begin{eqnarray}
\label{Hadmardadd}
\psi(x_t)-\psi(x)= \mathrm{d}\psi^*_t
(x_t-x)
\end{eqnarray}
holds.

For a smooth map $x_t$ with $x_t \rightarrow x$ as $t \downarrow0$,
the Hadamard differentiability of the function $\psi$ and the
additional condition (\ref{Hadmardadd}) imply that
%
\begin{eqnarray}
\label{Hadmardadd2}
\mathrm{d}\psi^*_t h \rightarrow \mathrm{d}\psi(x)h\qquad   \mbox{as }  t
\downarrow 0,
\end{eqnarray}
where the limit $\mathrm{d}\psi(x)$ is the Hadamard derivative of $\psi$ at $x$.

\textit{Note on additional condition}.
In many statistics applications, we
have the additional condition. 
For example, for functions $F(x)$ and $g(x)$, the map $\psi\dvtx F
\rightarrow\int g(x)\,\mathrm{d}F(x)$ satisfies the additional condition:
\[
\psi(F_t)-\psi(F)= \int g(x)\,\mathrm{d}(F_t-F) (x)
\]
here the map $\mathrm{d}\psi^*_t$ in (\ref{Hadmardadd}) is $\mathrm{d}\psi^*h= \int
g(x)\,\mathrm{d}h(x)$ which coincides with the Hadamard derivative of $\psi$.
For another example, consider a map $\psi\dvtx g \rightarrow(\int
g(x)\,\mathrm{d}F(x))^{-1}$. Then
\[
\psi(g_t)-\psi(g)= \frac{1}{\int g_t(x)\,\mathrm{d}F(x)}-\frac{1}{\int
g(x)\,\mathrm{d}F(x)}=
\frac{-\int[g_t(x)-g(x)]\,\mathrm{d}F(x)}{\int g_t(x)\,\mathrm{d}F(x)\int g(x)\,\mathrm{d}F(x)},
\]
and it shows the map $\psi$ satisfies the additional condition with
\[
\mathrm{d}\psi^*_t h = \frac{-\int h(x)\,\mathrm{d}F(x)}{\int g_t(x)\,\mathrm{d}F(x)\int g(x)\,\mathrm{d}F(x)}.
\]
If $g_t \rightarrow g$ as $t \downarrow0$, then
$\mathrm{d}\psi^*_t h$ converges to the Hadamard derivative of $\psi$:
\[
\mathrm{d}\psi h = \frac{-\int h(x)\,\mathrm{d}F(x)}{(\int g(x)\,\mathrm{d}F(x))^2}.
\]

\textit{Note on norm used in Theorem~\ref{mainthm} (below)}.
We treat the
set of c.d.f. functions ${\cal F}$ on ${\cal X}$ as a~subset of $\ell
^{\infty}({\cal X})$, the collection of all bounded functions on ${\cal X}$.
This means the norm on ${\cal F}$ is the sup-norm: for $F \in{\cal
F}$, $\|F\|=\sup_{x \in{\cal X}}|F(x)|$.
The convex subset $H$ of a Banach space ${\cal B}$ has the natural norm
from the Banach space and it is also denoted by $\|h\|$ for $h \in H$.
For all derivatives in the theorem, we use the operator norm.
The open subset $\Theta$ of $R^d$ has the Euclidean norm.

\begin{thm}\label{mainthm}
Suppose the map $(\theta,F,\eta) \rightarrow\Psi_{\theta,F}(\eta) \in
H$, $(\theta,F,\eta) \in\Theta\times{\cal F} \times H$, is:
\begin{enumerate}[(A3)]
\item[(A1)] Two times continuously differentiable with respect to
$\theta$ and two times Hadamard differentiable with respect to $\eta$
and Hadamard differentiable with respect to $F$ so that
the derivatives $\dot{\Psi}_{\theta,F}(\eta)$, $\ddot{\Psi}_{\theta
,F}(\eta)$,
$\mathrm{d}_{\eta}\Psi_{\theta,F}(\eta)$, $\mathrm{d}^2_{\eta}\Psi_{\theta,F}(\eta)$,
$\mathrm{d}_{\eta}\dot{\Psi}_{\theta,F}(\eta)$
and $\mathrm{d}_{F}\Psi_{\theta,F}(\eta)$ exist in some neighborhood of the true
value $(\theta_0,\eta_0,F_0)$ (where, e.g., $\dot{\Psi}_{\theta,F}(\eta
)$ is the first derivative with respect to $\theta$, and $\mathrm{d}_{\eta}\Psi
_{\theta,F}(\eta)$ is the first derivative with respect to $\eta$ in
the sense of Hadamard. Similarly, the rest is defined).
For each derivative, we assume the corresponding additional condition
(\ref{Hadmardadd}).

\item[(A2)] The true value $(\theta_0,\eta_0,F_0)$ satisfy $\eta_0 =
\Psi_{\theta_0,F_0}(\eta_0)$.
\item[(A3)] The linear operator $\mathrm{d}_{\eta}\Psi_{\theta_0,F_0}(\eta
_0)\dvtx {\cal B} \rightarrow{\cal B}$ has the operator norm $\|\mathrm{d}_{\eta}\Psi
_{\theta_0,F_0}(\eta_0)\|<1$.
\end{enumerate}
Then the solution $\eta_{\theta,F}$ to the equation
%
\begin{equation}
\label{ImplicitEquation2}
\eta= \Psi_{\theta,F}(\eta)
\end{equation}
exists in an neighborhood of $(\theta_0,F_0)$ and it is two times
continuously differentiable with respect to $\theta$ and Hadamard
differentiable with respect to $F$ in the neighborhood.
Moreover, the derivatives are given by
%
\begin{eqnarray}
\label{doteta}
\dot{\eta}_{\theta,F} &=& \bigl[I-\mathrm{d}_{\eta}
\Psi_{\theta,F}(\eta_{\theta
,F})\bigr]^{-1} \dot{
\Psi}_{\theta,F}(\eta_{\theta,F}),
\\
 \ddot{\eta}_{\theta,F} & = & \bigl[I-\mathrm{d}_{\eta}
\Psi_{\theta,F}(\eta_{\theta,F})\bigr]^{-1} \bigl[\ddot{\Psi
}_{\theta,F}(\eta_{\theta,F}) + \mathrm{d}_{\eta} \dot{
\Psi}_{\theta,F}(\eta _{\theta,F})\dot{\eta}_{\theta,F}^T
\nonumber
\\[-8pt]
\label{ddoteta}
\\[-8pt]
\nonumber
&&\hspace*{93pt}{}+ \mathrm{d}_{\eta} \dot{\Psi}_{\theta,F}^T(
\eta_{\theta
,F})\dot{\eta}_{\theta,F} + \mathrm{d}^2_{\eta}
\Psi_{\theta,F}(\eta_{\theta,F})\dot{\eta}_{\theta,F}\dot {
\eta}_{\theta,F}^T \bigr]
\end{eqnarray}
and
%
\begin{equation}
\label{dFeta}
\mathrm{d}_F\eta_{\theta,F}= \bigl[I-\mathrm{d}_{\eta}
\Psi_{\theta,F}(\eta_{\theta
,F})\bigr]^{-1}\,\mathrm{d}_F
\Psi_{\theta,F}(\eta_{\theta,F}).
\end{equation}
\end{thm}

\subsection{Proof of Theorem~\texorpdfstring{\protect\ref{mainthm}}{3.1}}\label{sec31}

We assumed the derivative $\mathrm{d}_{\eta}\Psi_{\theta_0,F_0}(\eta_0)$ exists
and its operator norm satisfies\break $\|\mathrm{d}_{\eta}\Psi_{\theta_0,F_0}(\eta_0)\|<1$.
By continuity of the map $(\theta,\eta,F) \rightarrow \mathrm{d}_{\eta}\Psi
_{\theta,F}(\eta)$, there are $\varepsilon>0$
and a neighborhood of $(\theta_0,\eta_0,F_0)$ such that
%
\begin{equation}
\label{detaPsi<1}
\bigl\|\mathrm{d}_{\eta}\Psi_{\theta,F}(\eta)\bigr\|<1-\varepsilon
\end{equation}
for all $(\theta,\eta,F)$ in the neighborhood.
In the following, we assume the parameters $(\theta,\eta,F)$ stay in
the neighborhood so that the inequality (\ref{detaPsi<1}) holds.

\textit{Existence and invertibility}.
Let $I\dvtx {\cal B}\rightarrow{\cal B}$ be the identity operator on the
space ${\cal B}$.
In the neighborhood discussed above, the map $(I-\mathrm{d}_{\eta} \Psi_{\theta
,F}(\eta))\dvtx{\cal B}\rightarrow{\cal B}$ has
the inverse $(I-\mathrm{d}_{\eta} \Psi_{\theta,F}(\eta))^{-1}$,
which is also a bounded linear map (cf. Kolmogorov and Fomin \cite{KF75}, Theorem~4, page 231).
It also follows that there is a neighborhood of $(\theta_0,\eta_0,F_0)$
such that, for each $(\theta,F)$,
the map $\eta\rightarrow\Psi_{\theta,F}(\eta)$ is a contraction
mapping in the neighborhood.
By Banach's contraction principle (cf. Agarwal, O'Regan and Sahu \cite{AOS09}, Theorem~4.1.5, page~178),
the solution to the equation (\ref{ImplicitEquation2}) exists uniquely
in the neighborhood.

\textit{Differentiability with respect to $F$}.
Fix $h$ 
in an appropriate space and let $F_t$
be a map such that
$F_{t=0}=F$,
$t^{-1}\{F_t-F\} \rightarrow h$
as $t \downarrow0$.
Then, $F_t \rightarrow F$
(as $t \downarrow0$).
We aim to find the limit of $t^{-1}\{\eta_{\theta,F_t} - \eta_{\theta
,F}\}$
as $t \downarrow0$.

(Step 1) First step is to show $\eta_{\theta,F_t} \rightarrow\eta
_{\theta,F}$ as $t \downarrow0$.
Due to equation (\ref{ImplicitEquation2}),
$\eta_{\theta,F}=\Psi_{\theta,F}(\eta_{\theta,F})$
and
$\eta_{\theta,F_t}=\Psi_{\theta,F_t}(\eta_{\theta,F_t})$.
It follows that\vspace*{-2pt}
%
\begin{eqnarray}
\{\eta_{\theta,F_t} - \eta_{\theta,F}\} & = & \bigl\{
\Psi_{\theta,F_t}(\eta_{\theta,F_t})-\Psi_{\theta,F}(\eta _{\theta,F})
\bigr\}
\nonumber
\\[-9pt]
\label{etadiff}\\[-9pt]
\nonumber
& = & \bigl\{\Psi_{\theta,F_t}(\eta_{\theta,F_t})-\Psi_{\theta,F_t}(\eta
_{\theta,F})\bigr\} +\bigl\{\Psi_{\theta,F_t}(\eta_{\theta,F})-
\Psi_{\theta,F}(\eta_{\theta,F})\bigr\}. 
\end{eqnarray}

Since the map $F \rightarrow\Psi_{\theta,F}(\eta)$ is continuous and
$F_t \rightarrow F$ (as $t \downarrow0$),
the second term in the right-hand side\vspace*{-2pt} is
\[
\Psi_{\theta,F_t}(\eta_{\theta,F})-\Psi_{\theta,F}(
\eta_{\theta,F}) =\mathrm{o}(1) \qquad \mbox{as }  t \downarrow 0.
\]
By the generalized Taylors theorem for Banach spaces (cf. \cite{Zeidler95}, page 243, Theorem~4C), the first term in the right-hand
side\vspace*{-2pt} is
\begin{eqnarray*}
\bigl\|\Psi_{\theta,F_t}(\eta_{\theta,F_t})-\Psi_{\theta,F_t}(
\eta_{\theta
,F})\bigr\| & \leq& \sup_{\tau\in[0,1]} \bigl\|\mathrm{d}_\eta
\Psi_{\theta,F_t}\bigl(\eta_{\theta
,F}+\tau(\eta_{\theta,F_t}-
\eta_{\theta,F})\bigr)\bigr\| \|\eta_{\theta,F_t}-\eta _{\theta,F}\|
\\
& \leq& (1-\varepsilon) \|\eta_{\theta,F_t}-\eta_{\theta,F}\|,
\end{eqnarray*}
where the last inequality is due to (\ref{detaPsi<1}).

It follows from (\ref{etadiff})\vspace*{-2pt} that
\[
\|\eta_{\theta,F_t} - \eta_{\theta,F}\|  \leq \mathrm{o}(1) +(1-\varepsilon) \|
\eta_{\theta,F_t}-\eta_{\theta,F}\| \qquad  \mbox{as }  t \downarrow 0.
\]
This shows $\eta_{\theta,F_t} \rightarrow\eta_{\theta,F}$ as $t
\downarrow0$.

(Step 2)
By the Hadamard differentiability of the map $F \rightarrow\Psi_{\theta
,F}(\eta)$
and the additional condition ((\ref{Hadmardadd}) and (\ref{Hadmardadd2})),
there is a linear operator $h \rightarrow \mathrm{d}_F \Psi^*_th$ such that the
first term in the right-hand side of (\ref{etadiff}) can be expressed as
\[
\bigl\{\Psi_{\theta,F_t}(\eta_{\theta,F_t})-\Psi_{\theta,F}(
\eta_{\theta
,F_t})\bigr\} = \mathrm{d}_F \Psi^*_t
(F_t-F),
\]
and\vspace*{-3pt}
\[
\mathrm{d}_F \Psi^*_t \rightarrow \mathrm{d}_F
\Psi_{\theta,F}(\eta_{\theta,F}) \qquad \mbox{as } t \downarrow0.
\]

Similarly, there is a linear operator $h' \rightarrow \mathrm{d}_{\eta} \Psi
^*_th'$ such that the second term in the right-hand side of (\ref{etadiff})\vspace*{-2pt} is
\[
\bigl\{\Psi_{\theta,F}(\eta_{\theta,F_t})-\Psi_{\theta,F}(
\eta_{\theta,F})\bigr\} 
 =  \mathrm{d}_{\eta}
\Psi^*_t \{\eta_{\theta,F_t}-\eta_{\theta,F}\}
\]
and\vspace*{-3pt}
\[
\mathrm{d}_{\eta} \Psi^*_t \rightarrow \mathrm{d}_{\eta}
\Psi_{\theta,F}(\eta_{\theta
,F}) \qquad \mbox{as } t \downarrow0.
\]

Altogether, equation (\ref{etadiff}) can be written\vspace*{-2pt} as
\[
\{\eta_{\theta,F_t} - \eta_{\theta,F}\}  =  \mathrm{d}_F
\Psi^*_t(F_t-F) +\mathrm{d}_{\eta} \Psi^*_t
\{\eta_{\theta,F_t}-\eta_{\theta,F}\}.
\]
It follows that
\[
\bigl[I-\mathrm{d}_{\eta} \Psi^*_t\bigr]\{\eta_{\theta,F_t}-
\eta_{\theta,F}\}= \mathrm{d}_F \Psi^*_t(F_t-F),
\]
where $I$ is the identity operator in the space ${\cal B}$.

Since we have the inequality (\ref{detaPsi<1})
and $\mathrm{d}_{\eta} \Psi^*_t \rightarrow \mathrm{d}_{\eta} \Psi_{\theta,F}(\eta_{\theta
,F})$ as $t \downarrow0$,
the inverse $[I-\mathrm{d}_{\eta} \Psi^*_t]^{-1}$ exists for small $t> 0$.
Therefore, when $t^{-1}(F_t-F) \rightarrow h$ as $t \downarrow0$, we
have that
\begin{eqnarray*}
t^{-1}\{\eta_{\theta,F_t}-\eta_{\theta,F}\} & = &
\bigl[I-\mathrm{d}_{\eta} \Psi^*_t\bigr]^{-1}\,
\mathrm{d}_F \Psi^*_tt^{-1}(F_t-F)
\\
& \rightarrow& \bigl[I-\mathrm{d}_{\eta} \Psi_{\theta,F}(\eta_{\theta,F})
\bigr]^{-1}\,\mathrm{d}_F \Psi_{\theta,F}(\eta_{\theta,F})h \qquad  \mbox{as }  t \downarrow0.
\end{eqnarray*}

Since the limit is a bounded and linear map of $h$,
the function $\eta_{\theta,F}(x)$ is Hadamard differentiable with
respect to $F$ with the derivative
\[
\mathrm{d}_F\eta_{\theta,F}=\bigl[I-\mathrm{d}_{\eta}
\Psi_{\theta,F}(\eta_{\theta
,F})\bigr]^{-1}\,\mathrm{d}_F
\Psi_{\theta,F}(\eta_{\theta,F}).
\]

\textit{Differentiability with respect to $\theta$}.
Similar proof as above can show that, for $t^{-1}(\theta_t-\theta)
\rightarrow a \in R^d$ as $t \downarrow0$, we have
\[
t^{-1}\{\eta_{\theta_t,F} - \eta_{\theta,F}\}
 \rightarrow \bigl[I-\mathrm{d}_{\eta}
\Psi_{\theta,F}(\eta_{\theta,F})\bigr]^{-1} a^T\dot{
\Psi}_{\theta,F}(\eta_{\theta,F}).
\]
It follows that the first derivative $\dot{\eta}_{\theta,F}$ of $\eta
_{\theta,F}(x)$ with respect to $\theta$ is given by
%
\begin{equation}
\label{dotg1}
a^T\dot{\eta}_{\theta,F}=\bigl[I-\mathrm{d}_{\eta}
\Psi_{\theta,F}(\eta_{\theta
,F})\bigr]^{-1} a^T\dot{
\Psi}_{\theta,F}(\eta_{\theta,F}).
\end{equation}

Now we show the second derivative  of $\eta_{\theta,F}(x)$ with
respect to $\theta$.
From (\ref{dotg1}), we have
\[
a^T\dot{\eta}_{\theta,F}= a^T\dot{
\Psi}_{\theta,F}(\eta_{\theta
,F})+\mathrm{d}_{\eta} \Psi_{\theta,F}(
\eta_{\theta,F}) \bigl(a^T\dot{\eta}_{\theta,F}\bigr).
\]

Using this equation, for $t^{-1}(\theta_t-\theta) \rightarrow b \in
R^d$ as $t \downarrow0$,
\begin{eqnarray*}
& & t^{-1}\bigl\{a^T\dot{\eta}_{\theta_t,F} -
a^T\dot{\eta}_{\theta,F}\bigr\}
\\
&& \quad  =  t^{-1}\bigl\{a^T\dot{\Psi}_{\theta_t,F}(
\eta_{\theta_t,F})-a^T\dot{\Psi }_{\theta,F}(
\eta_{\theta,F})\bigr\} \\
&& \qquad {}+ t^{-1}\bigl\{\mathrm{d}_{\eta}
\Psi_{\theta
_t,F}(\eta_{\theta_t,F}) \bigl(a^T\dot{
\eta}_{\theta_t,F}\bigr)-\mathrm{d}_{\eta} \Psi _{\theta,F}(
\eta_{\theta,F}) \bigl(a^T\dot{\eta}_{\theta,F}\bigr)\bigr\}
\\
&&\quad  =  t^{-1}\bigl\{a^T\dot{\Psi}_{\theta_t,F}(
\eta_{\theta_t,F})-a^T\dot{\Psi }_{\theta,F}(
\eta_{\theta_t,F})\bigr\} + t^{-1}\bigl\{a^T\dot{
\Psi}_{\theta
,F}(\eta_{\theta_t,F})-a^T\dot{
\Psi}_{\theta,F}(\eta_{\theta,F})\bigr\}
\\
&&\qquad {}+ t^{-1}\bigl\{\mathrm{d}_{\eta} \Psi_{\theta_t,F}(
\eta_{\theta_t,F}) \bigl(a^T\dot{\eta }_{\theta_t,F}
\bigr)-\mathrm{d}_{\eta} \Psi_{\theta,F}(\eta_{\theta_t,F})
\bigl(a^T\dot {\eta}_{\theta_t,F}\bigr)\bigr\}
\\
&&\qquad  {}+ t^{-1}\bigl\{\mathrm{d}_{\eta} \Psi_{\theta,F}(
\eta_{\theta_t,F}) \bigl(a^T\dot{\eta }_{\theta_t,F}
\bigr)-\mathrm{d}_{\eta} \Psi_{\theta,F}(\eta_{\theta,F})
\bigl(a^T\dot{\eta }_{\theta_t,F}\bigr)\bigr\}
\\
& & \qquad {}+ t^{-1}\bigl\{\mathrm{d}_{\eta} \Psi_{\theta,F}(
\eta_{\theta,F}) \bigl(a^T\dot{\eta }_{\theta_t,F}
\bigr)-\mathrm{d}_{\eta} \Psi_{\theta,F}(\eta_{\theta,F})
\bigl(a^T\dot{\eta }_{\theta,F}\bigr)\bigr\}. 
\end{eqnarray*}

By the differentiability with respect to $\theta$,
the each term in the right-hand side has the limit as follows, as $t
\downarrow0$,
\begin{eqnarray*}
t^{-1}\bigl\{a^T\dot{\Psi}_{\theta_t,F}(
\eta_{\theta_t,F})-a^T\dot{\Psi }_{\theta,F}(
\eta_{\theta_t,F})\bigr\} & \rightarrow&a^T\ddot{
\Psi}_{\theta,F}(\eta_{\theta,F})b,
\\
t^{-1}\bigl\{a^T\dot{\Psi}_{\theta,F}(
\eta_{\theta_t,F})-a^T\dot{\Psi }_{\theta,F}(
\eta_{\theta,F})\bigr\} & \rightarrow& a^T \,\mathrm{d}_{\eta} \dot{
\Psi}_{\theta,F}(\eta_{\theta,F}) \bigl(\dot {\eta}_{\theta,F}^Tb
\bigr),
\\
t^{-1}\bigl\{\mathrm{d}_{\eta} \Psi_{\theta_t,F}(\eta_{\theta_t,F})
\bigl(a^T\dot{\eta }_{\theta_t,F}\bigr)-\mathrm{d}_{\eta}
\Psi_{\theta,F}(\eta_{\theta_t,F}) \bigl(a^T\dot {
\eta}_{\theta_t,F}\bigr)\bigr\} & \rightarrow& \bigl\{\mathrm{d}_{\eta} \dot{
\Psi}_{\theta,F}(\eta_{\theta
,F}) \bigl(a^T\dot{
\eta}_{\theta,F}\bigr)\bigr\}^Tb,
\\
t^{-1}\bigl\{\mathrm{d}_{\eta} \Psi_{\theta,F}(\eta_{\theta_t,F})
\bigl(a^T\dot{\eta }_{\theta_t,F}\bigr)-\mathrm{d}_{\eta}
\Psi_{\theta,F}(\eta_{\theta,F}) \bigl(a^T\dot{\eta
}_{\theta_t,F}\bigr)\bigr\} & \rightarrow& \mathrm{d}^2_{\eta}
\Psi_{\theta,F}(\eta_{\theta,F}) \bigl(a^T\dot{\eta
}_{\theta,F}\bigr) \bigl(\dot{\eta}_{\theta,F}^Tb\bigr),\\[-27pt]
\end{eqnarray*}
\begin{eqnarray*}
 && t^{-1}\bigl\{\mathrm{d}_{\eta} \Psi_{\theta,F}(\eta_{\theta,F})
\bigl(a^T\dot{\eta }_{\theta_t,F}\bigr)-\mathrm{d}_{\eta}
\Psi_{\theta,F}(\eta_{\theta,F}) \bigl(a^T\dot{\eta
}_{\theta,F}\bigr)\bigr\}
\\
 && \quad =  \mathrm{d}_{\eta} \Psi_{\theta,F}(
\eta_{\theta,F}) t^{-1}\bigl\{a^T\dot{\eta
}_{\theta_t,F}- a^T\dot{\eta}_{\theta,F}\bigr\},
\end{eqnarray*}
where the last equality is due to the linearity of the operator $\mathrm{d}_{\eta
} \Psi_{\theta,F}(\eta_{\theta,F})\dvtx{\cal B}\rightarrow{\cal B}$
(the Hadamard derivative of $\Psi_{\theta,F}(\eta_{\theta,F})$ with
respect to $\eta$).

Using additional condition and the Hadamard differentiability in (A1),
by similar argument to the case for the differentiability with respect
to $F$,
we can show that
\begin{eqnarray*}
& & t^{-1}\bigl\{a^T\dot{\eta}_{\theta_t,F} -
a^T\dot{\eta}_{\theta,F}\bigr\}
\\
&& \quad =  a^T\ddot{\Psi}_{\theta,F}(\eta_{\theta,F})b
+a^T \,\mathrm{d}_{\eta} \dot {\Psi}_{\theta,F}(\eta_{\theta,F})
\bigl(\dot{\eta}_{\theta,F}^Tb\bigr) + \bigl\{\mathrm{d}_{\eta}
\dot{\Psi}_{\theta,F}(\eta_{\theta,F}) \bigl(a^T\dot{\eta
}_{\theta,F}\bigr)\bigr\}^Tb
\\
&& \qquad {}+ \mathrm{d}^2_{\eta} \Psi_{\theta,F}(\eta_{\theta,F})
\bigl(a^T\dot{\eta}_{\theta
,F}\bigr) \bigl(\dot{
\eta}_{\theta,F}^Tb\bigr) + \mathrm{d}_{\eta} \Psi_{\theta,F}(
\eta_{\theta,F}) t^{-1}\bigl\{a^T\dot{\eta
}_{\theta_t,F}- a^T\dot{\eta}_{\theta,F}\bigr\}+\mathrm{o}(1).
\end{eqnarray*}

By rearranging this, we obtain
\begin{eqnarray*}
& & \bigl[I-\mathrm{d}_{\eta} \Psi_{\theta,F}(\eta_{\theta,F})
\bigr]t^{-1}\bigl\{a^T\dot{\eta }_{\theta_t,F} -
a^T\dot{\eta}_{\theta,F}\bigr\}
\\
&&\quad  =  a^T\ddot{\Psi}_{\theta,F}(\eta_{\theta,F})b
+a^T\, \mathrm{d}_{\eta} \dot {\Psi}_{\theta,F}(\eta_{\theta,F})
\bigl(\dot{\eta}_{\theta,F}^Tb\bigr) + \bigl\{\mathrm{d}_{\eta}
\dot{\Psi}_{\theta,F}(\eta_{\theta,F}) \bigl(a^T\dot{\eta
}_{\theta,F}\bigr)\bigr\}^Tb
\\
& &\qquad  {}+ \mathrm{d}^2_{\eta} \Psi_{\theta,F}(\eta_{\theta,F})
\bigl(a^T\dot{\eta}_{\theta
,F}\bigr) \bigl(\dot{
\eta}_{\theta,F}^Tb\bigr) +\mathrm{o}(1),
\end{eqnarray*}
and hence, as $t \downarrow0$,
\[
t^{-1}\bigl\{a^T\dot{\eta}_{\theta_t,F} - a^T
\dot{\eta}_{\theta,F}\bigr\} \rightarrow a^T\ddot{
\eta}_{\theta,F}b,
\]
where
\begin{eqnarray*}
a^T\ddot{\eta}_{\theta,F}b & = & \bigl[I-\mathrm{d}_{\eta}
\Psi_{\theta,F}(\eta_{\theta,F})\bigr]^{-1} \bigl[
a^T\ddot {\Psi}_{\theta,F}(\eta_{\theta,F})b +a^T\,
\mathrm{d}_{\eta} \dot{\Psi}_{\theta
,F}(\eta_{\theta,F}) \bigl(\dot{
\eta}_{\theta,F}^Tb\bigr)
\\
& & {}+ \bigl\{\mathrm{d}_{\eta} \dot{\Psi}_{\theta,F}(
\eta_{\theta
,F}) \bigl(a^T\dot{\eta}_{\theta,F}\bigr)\bigr
\}^Tb + \mathrm{d}^2_{\eta} \Psi_{\theta,F}(
\eta_{\theta,F}) \bigl(a^T\dot{\eta}_{\theta
,F}\bigr) \bigl(
\dot{\eta}_{\theta,F}^Tb\bigr) \bigr].
\end{eqnarray*}

Therefore, $\dot{\eta}_{\theta,F}$ is differentiable with respect to
$\theta$ with derivative $\ddot{\eta}_{\theta,F}$.\\


\section{Example~1 continued}\label{sec4}

As an application of the main result (Theorem~\ref{mainthm}), we show
existence and differentiability of solution to the operator equation in
Example~1.

\begin{thm}
Suppose that
%
\begin{equation}
\label{ConditionA3}
E_F \biggl(\frac{\delta}{1+\delta}W^2(s;
\beta,A) \biggr)> \operatorname{Var}_FW(s;\beta,A),
\end{equation}
where $\operatorname{Var}_FW(s;\beta,A) = E_FW^2(s;\beta,A)-\{E_FW(s;\beta,A)\}^2$.
Then the solution $A_{\beta,F}(t)$ to the operator equation
\[
A=\Psi_{\beta,F}(A)
\]
exists in an neighborhood of $(\beta_0,F_0)$
and it is two times continuously differentiable with respect to $\beta$
and Hadamard differentiable with respect to $F$ in the neighborhood,
where the operator $\Psi_{\beta,F}(A)$ is given in (\ref{PsiPodds}).
\end{thm}

For the proof, we verify conditions (A1), (A2) and (A3) in Theorem~\ref{mainthm} so that the differentiability of the solution
is implied by the theorem.

\textit{Verification of condition (A1).}
We show that the map $\Psi_{\beta,F}(A)$ defined by (\ref{PsiPodds}) is
differentiable with respect to $\beta$, $F$ and $A$.

(\textit{The derivative of $\Psi_{\beta,F}(A)$ with respect to $F$})
Suppose a map $t \rightarrow F_t$ satisfies
$t^{-1}(F_t-F) \rightarrow h$ as $t \downarrow0$.
\begin{eqnarray*}
t^{-1}\bigl\{\Psi_{\beta,F_t}(A)- \Psi_{\beta,F}(A)\bigr\}
& = & t^{-1} \biggl\{ E_{F_t}\int
_0^u \frac{\mathrm{d}N(s)}{E_{F_t} W(s;\beta,A)} - E_F\int
_0^u \frac{\mathrm{d}N(s)}{E_F W(s;\beta,A)} \biggr\}
\\
& = & t^{-1} \biggl\{ E_{F_t}\int_0^u
\frac{\mathrm{d}N(s)}{E_{F_t} W(s;\beta,A)}-E_{F}\int_0^u
\frac
{\mathrm{d}N(s)}{E_{F_t} W(s;\beta,A)} \biggr\}
\\
&& {}+t^{-1} \biggl\{ E_{F}\int_0^u
\frac{\mathrm{d}N(s)}{E_{F_t} W(s;\beta,A)} - E_F\int_0^u
\frac{\mathrm{d}N(s)}{E_F W(s;\beta,A)} \biggr\}.
\end{eqnarray*}
After a simple calculation the right-hand side is equal to
\begin{eqnarray}
&& \mathrm{d}
\Psi_t^* \bigl(t^{-1}\{F_t-F\}\bigr)
\nonumber
\\[-8pt]
\label{example1additional}\\[-8pt]
\nonumber
&& \quad  =  E_{t^{-1}\{F_t-F\}}\int_0^u
\frac{\mathrm{d}N(s)}{E_{F_t} W(s;\beta,A)} - E_{F}\int_0^u
\frac{ E_{t^{-1}\{F_t-F\}} W(s;\beta,A) }{E_F W(s;\beta
,A)E_{F_t} W(s;\beta,A)} \,\mathrm{d}N(s),
\nonumber
\end{eqnarray}
where the notation $E_F f$ means $\int f \,\mathrm{d}F$.
The expression (\ref{example1additional}) shows the additional
condition~(\ref{Hadmardadd}) is satisfied.
Moreover, as $t \downarrow0$, the expression converges to
\begin{eqnarray*}
\mathrm{d}_F\Psi_{\beta,F}(A)h=E_{h}\int
_0^u \frac{\mathrm{d}N(s)}{E_{F} W(s;\beta,A)} - E_{F}\int
_0^u \frac{E_{h} W(s;\beta,A) }{\{E_F W(s;\beta,A)\}^2}\,\mathrm{d}N(s).
\end{eqnarray*}
This shows the map $F \rightarrow\Psi_{\beta,F}(A)$ is Hadamard
differentiable at $(\beta,A,F)$
with derivative $\mathrm{d}_{F}\Psi_{\beta,F}(A)$ and additional condition satisfied
(clearly, the derivative is linear in $h$, we omit the proof of
boundedness of $\mathrm{d}_{F}\Psi_{\beta,F}(A)$).


For the rest the derivatives, the proofs are similar and
straightforward, therefore, we omit the proof and just give the
derivatives in Appendix~\ref{appB}.

\textit{Verification of condition (A2)}.
Let $F_0$ be the true c.d.f. and $\beta_0$ be the true value of $\beta$.
Since the true value $A_0$ of $A$ is the maximizer of the expected
log-likelihood
\[
\int\bigl\{\delta\bigl(\beta_0'Z + \log a(U)\bigr)-(1+
\delta)\log\bigl(1+\mathrm{e}^{\beta
_0'Z}A(U)\bigr)\bigr\} \,\mathrm{d}F_0,
\]
the same method to derive the equation (\ref{hatAbeta}) can be applied
to show
\[
A_0(u)=\int_0^u
\frac{E_{F_0} \,\mathrm{d}N(s) }{E_{F_0} W(s;\beta_0,A_0)}=\Psi _{\beta_0,F_0}(A_0),
\]
where
$E_{F_0}\,\mathrm{d}N(s)= \int \mathrm{d}N(s) \,\mathrm{d}F_0$, $E_{F_0} W(s;\beta_0,A_0)=\int
W(s;\beta_0,A_0)\,\mathrm{d}F_0$
and $\Psi_{\beta,F}(A)$ is defined in (\ref{PsiPodds}).

\textit{Verification of condition (A3)}.
The derivatives $\mathrm{d}_A\Psi_{\beta,F}(A)$ and $\mathrm{d}_A W(s;\beta,A)$ are given
in (\ref{dAPsi}) and (\ref{dAW}), respectively, in Appendix~\ref{appB}.
We consider the $sup$-norm on the space of total variation bounded
cadlag functions $h_1(u)$ on $[0,\tau]$.
For all $h_1(u)$ such that $\|h_1(u)\|=\sup_{u\in[0,\tau]} |h_1(u)|
\leq1$, 
we have that
\begin{eqnarray*}
\bigl|\mathrm{d}_A W(s;\beta,A)h_1\bigr| & \leq& \frac{(1+\delta)\mathrm{e}^{2\beta'Z}Y(s) |h_1(U)|}{\{1+\mathrm{e}^{\beta
'Z}A(U)\}^2} \leq
\frac{(1+\delta)\mathrm{e}^{2\beta'Z}Y(s)}{\{1+\mathrm{e}^{\beta'Z}A(U)\}^2}
\\
& \leq& \frac{(1+\delta)^2\mathrm{e}^{2\beta'Z}Y(s)}{\{1+\mathrm{e}^{\beta'Z}A(U)\}^2} = W^2(s;\beta,A).
\end{eqnarray*}
We assumed $P\{\delta Y(s)=1\}>0$ for each $s \in[0,\tau]$
so that the last inequality in the above equation is strict inequality
with positive probability for each $s$.
This implies
%
\begin{equation}
\label{W^2-dAW} W^2(s;\beta,A) - \bigl|\mathrm{d}_A W(s;
\beta,A)h_1\bigr| \geq \frac{\delta}{1+\delta
}W^2(s;\beta,A)>0
\end{equation}
with positive probability for each $s$.

Then, by (\ref{ConditionA3}) and (\ref{W^2-dAW}), we have that, for
each $s$,
\[
E_FW^2(s;\beta,A) - E_F \bigl|\mathrm{d}_A
W(s;\beta,A)h_1\bigr| > E_FW^2(s;\beta,A)-\bigl\{ E_FW(s;\beta,A)\bigr\}^2>0.
\]
It follows that
\[
\bigl|\mathrm{d}_A\Psi_{\beta,F}(A)h_1 \bigr|  \leq
E_{F}\int_0^u \frac{E_{F} |\mathrm{d}_AW(s;\beta,A)h_1| }{\{E_{F}
W(s;\beta,A)\}^2 }\,\mathrm{d}N(s)
< E_{F}\int_0^u
\,\mathrm{d}N(s) \leq1.
\]
This demonstrates the operator $h_1 \rightarrow \mathrm{d}_A W(s;\beta,A)h_1$
has the operator norm smaller than one.

We have completed verification of conditions (A1), (A2) and (A3) in
Theorem~\ref{mainthm}.
By the theorem it follows that the derivatives of the function (\ref
{hatAbeta}) is given by
equations (\ref{doteta}), (\ref{ddoteta}) and~(\ref{dFeta}) (needs
replacement $\theta$ with $\beta$ and $\eta$ with $A$).

\section{Example~2 continued}\label{sec5}

The generic\vspace*{1pt} form of c.d.f. for combined samples is
$F=\sum_{s=1}^2w_sF_s$
where $w_s>0$, $s=1,2$, and $w_1+w_2=1$ and $F_1, F_2$ are c.d.f.s for the
samples in $V$ and $\overline{V}$, respectively.

For $\theta\in R^d$, $F$ and function $g(x)$, define
%
\begin{equation}
\label{PsithetaF} \Psi_{\theta,F}(g) 
=
\frac{\partial_x \int \pi_1(\mathrm{d}F)}{A(x;\theta,g,F)},
\end{equation}
where $\pi_s\dvtx F=\sum_{s'=1}^2w_{s'}F_{s'} \rightarrow w_sF_s$, $s=1,2$,
are projections,
and
%
\begin{equation}
\label{A} A(x;\theta,g,F)=1 - \int \frac{ f(y|x;\theta)}{f_Y(y;\theta,g)} \pi_2
(\mathrm{d}F).
\end{equation}

Then the function $g_{\theta,F_n}(x)$ given by (\ref{hatg}) is the
solution to the operator equation
%
\begin{equation}
\label{g=Psig} g(x)=\Psi_{\theta,F}(g) (x)
\end{equation}
with $F=F_n$.

We show the differentiability of the solution $g_{\theta,F}(x)$ to the
equation (\ref{g=Psig}) with respect to $\theta$ and $F$.

\begin{thm}
Let $\theta_0$, $g_0$ and $F_0=\sum_{s=1}^2 w_{s0}F_{s0}$ be the true
values of $\theta$, $g$ and $F$ at which data are generated.
We assume that
%
\begin{equation}
\label{ConditionA3ex2} \frac{w_{20}}{w_{10}}<1
\end{equation}
and the function $f(y|x;\theta)$ is twice continuously differentiable
with respect to $\theta$.
Then the solution $g_{\theta,F}(x)$ to the operator equation (\ref{g=Psig}) exists in an neighborhood of $(\theta_0,F_0)$ and it is two
times continuously differentiable with respect to $\theta$ and Hadamard
differentiable with respect to $F$ in the neighborhood.
\end{thm}

To prove the theorem, we verify conditions (A1), (A2) and (A3) in
Theorem~\ref{mainthm} so that the results follows from that theorem.

We denote $f=f(y|x;\theta)$, $f_Y=f_Y(y;\theta,g)$, $A=A(x;\theta,g,F)$,
$\dot{f}=\frac{\partial}{\partial\theta}f(y|x;\theta)$,
$\ddot{f}=\frac{\partial^2}{\partial\theta\, \partial\theta
^T}f(y|x;\theta)$,
$\dot{f}_Y=\int\dot{f}(y|x;\theta)g(x)\,\mathrm{d}x$,\vspace*{1pt} and $\ddot{f}_Y=\int\ddot
{f}(y|x;\theta)g(x)\,\mathrm{d}x$.

\textit{Verification of condition (A1)}.
We show that the map $\Psi_{\theta,F}(g)$ is differentiable with
respect to $\theta$, $F$ and $g$.

(\textit{The derivative of $\Psi_{\theta,F}(g)$ with respect to $F$})
Suppose a map $t \rightarrow F_t$ satisfies
$t^{-1}(F_t-F) \rightarrow h$ as $t \downarrow0$.

Then
\begin{eqnarray*}
& & \Psi_{\theta,F_t}(g)- \Psi_{\theta,F}(g) \\
&& \quad = \frac{\partial_x\int\pi_1(\mathrm{d}F_t)}{A(x;\theta,g,F_t)} -
\frac{\partial_x \int\pi_1(\mathrm{d}F)}{A(x;\theta,g,F)}
\\
& & \quad  = \frac{
 (\partial_x\int\pi_1 [\mathrm{d} (F_t-F)] ) A(x;\theta,g,F)
- (\partial_x \int\pi_1(\mathrm{d}F) ) \{A(x;\theta,g,F_t)-A(x;\theta
,g,F)\}
}{A(x;\theta,g,F_t)A(x;\theta,g,F)}.
\end{eqnarray*}

By equation (\ref{A}), the right-hand side is equal to
\begin{eqnarray*}
&&\hspace*{-4pt} \mathrm{d}_F\Psi^*_t(g) (F_t-F) \\
&& \hspace*{-7.5pt}\quad =
\frac{  (\partial_x\int\pi_1 [\mathrm{d} (F_t-F)] )A(x;\theta
,g,F)+ (\partial_x\int\pi_1(\mathrm{d}F) )\int { f(y|x;\theta
)}/{f_Y(y;\theta,g)} \pi_2 [\mathrm{d} (F_t-F)]}{
A(x;\theta,g,F_t)A(x;\theta,g,F)}.
\end{eqnarray*}
This shows the additional condition (\ref{Hadmardadd}) is satisfied.
Moreover, as $t \downarrow0$,
\[
t^{-1} \bigl\{ \Psi_{\theta,F_t}(g)- \Psi_{\theta,F}(g) \bigr\}
=  t^{-1}\,\mathrm{d}_F\Psi^*_t(g) (F_t-F)
\rightarrow \mathrm{d}_F\Psi_{\theta,F}(g)h,
\]
where the map $\mathrm{d}_F\Psi_{\theta,F}(g)$ is given by
\[
\label{dFPsi}
\mathrm{d}_F\Psi_{\theta,F}(g)h=\frac{ (\partial_x \int\pi_1 (\mathrm{d}h)
)A(x;\theta,g,F)
+ (\partial_x \int\pi_1(\mathrm{d}F) )
\int { f(y|x;\theta)}/{f_Y(y;\theta,g)} \pi_2 (\mathrm{d}h)
}{\{A(x;\theta,g,F)\}^2}.
\]

Hence, the map $F \rightarrow\Psi_{\theta,F}(g)$ is Hadamard
differentiable at $(\theta,g,F)$ with derivative $\mathrm{d}_{F}\Psi_{\theta,F}(g)$
(clearly, the derivative is linear in $h$, we omit the proof of
boundedness of $\mathrm{d}_{F}\Psi_{\theta,F}(g)$).

Similarly, other (Hadamard) differentiability of map can be shown. In
Appendix~\ref{appC}, we list the derivatives without proofs.

\textit{Verification of condition (A2)}.
To verify (A2), we show that, at $(\theta_0,F_0)$, $g_0(x)$ is a
solution to the operator equation (\ref{g=Psig}).

Since $\partial_x \int \mathrm{d}F_{10}=\int f(y|x;\theta_0)g_0(x)\,\mathrm{d}y=g_0(x)$,
and $\frac{\mathrm{d}F_{20}(y)}{\mathrm{d}y}=f_Y(y;\theta_0,g_0)$,
$w_{10}+w_{20}=1$, we have
%
\begin{eqnarray}
\Psi_{\theta_0,F_0}(g_0) (x)&=& \frac{ w_{10} \partial_x \int \mathrm{d}F_{10}}{1 -
w_{20} \int { f(y|x;\theta_0)}/{f_Y(y;\theta_0,g_0)} \,\mathrm{d}F_{20}}
\nonumber
\\[-8pt]
\label{VerifyA2ex2}\\[-8pt]
\nonumber
& =&
\frac{ w_{10} g_0(x)}{1 - w_{20} \int { f(y|x;\theta
_0)}/{f_Y(y;\theta_0,g_0)}f_Y(y;\theta_0,g_0)\,\mathrm{d}y}=g_0(x),
\end{eqnarray}
where we used $\int f(y|x;\theta)\,\mathrm{d}y=1$ for each $x$.

\textit{Verification of condition (A3).}
Let $L_1$ be the space of all real valued measurable functions $h(x)$
with $\|h\|_1=\int|h(x)|\,\mathrm{d}x < \infty$.
Then $L_1$ is a Banach space with the norm $\| \cdot\|_1$. The
sup-norm is denoted by $\|h\|_{\infty}=\sup_{x}|h(x)|$.

The derivatives $\mathrm{d}_g\Psi_{\theta,F}(g)$ and $\mathrm{d}_gA(x;\theta,g,F)$ are,
respectively, given in (\ref{dgPsi}) and (\ref{dgA}).

Since $\partial_x \int\pi_1(\mathrm{d}F_0)=w_{10} g_0(x)$, (\ref{dgPsi}) implies
\[
\mathrm{d}_g\Psi_{\theta_0,F_0}(g_0)h^*  =
\frac{- w_{10}g_0(x) \,\mathrm{d}_gA(x;\theta_0,g_0,F_0)h^*}{\{A(x;\theta
_0,g_0,F_0)\}^2}.
\]

By (\ref{A}) together with $\pi_2(\mathrm{d}F_0)=w_{20}f_Y(y;\theta_0,g_0)\,\mathrm{d}y$,
and $\int f(y|x;\theta)\,\mathrm{d}y=1$, for all $x$, we have
\[
A(x;\theta,g_0,F_0)  =  1 - \int \frac{ f(y|x;\theta_0)}{f_Y(y;\theta_0,g_0)}
\pi_2 (\mathrm{d}F_0) = 1 - w_{20}=w_{10}.
\]

These equations and (\ref{dgA}) imply
%
\begin{equation}
\label{dgPsi0}
\mathrm{d}_g\Psi_{\theta_0,F_0}(g_0)h^*  =  -
\frac{w_{20}}{w_{10}}g_0(x) \int f(y|x;\theta_0)
\frac{\int
f(y|x;\theta_0)h^*(x)\,\mathrm{d}x}{f_Y(y;\theta_0,g_0)}\, \mathrm{d}y.
\end{equation}
The $L_1$ norm of (\ref{dgPsi0}) is
\begin{eqnarray*}
\bigl\|\mathrm{d}_g\Psi_{\theta_0,F_0}(g_0)h^*\bigr\|_1 & =
& \int\biggl\llvert \frac{w_{20}}{w_{10}}g_0(x) \int f(y|x;
\theta_0) \frac
{\int f(y|x;\theta_0)h^*(x)\,\mathrm{d}x }{f_Y(y;\theta_0,g_0)} \,\mathrm{d}y\biggr\rrvert \,\mathrm{d}x
\\
& \leq& \frac{w_{20}}{w_{10}} \int g_0(x) \biggl(\int f(y|x;
\theta_0) \frac{\int f(y|x;\theta_0)|h^*(x)|\,\mathrm{d}x }{f_Y(y;\theta_0,g_0)} \,\mathrm{d}y \biggr) \,\mathrm{d}x
\\
& = & \frac{w_{20}}{w_{10}} \int\bigl|h^*(x)\bigr|\,\mathrm{d}x \qquad  \biggl(
\mbox{by Fubini's theorem and } \int f(y|x;\theta_0)\,\mathrm{d}y=1
\biggr)
\\
& = & \frac{w_{20}}{w_{10}}\|h^*\|_1.
\end{eqnarray*}
From the calculation above, we see that the operator $h^* \rightarrow
\mathrm{d}_g\Psi_{\theta_0,F_0}(g_0)h^*$ has the operator norm ${\leq}\frac
{w_{20}}{w_{10}}$.
Since we assumed $\frac{w_{20}}{w_{10}}<1$, we have condition (A3).

\section{Asymptotic normality of maximum profile likelihood~estimator}\label{sec6}

Hirose \cite{Hirose10} showed the efficiency of the maximum profile likelihood
estimator in semi-parametric models
using the direct asymptotic expansion of the profile likelihood.
The method gives alternative to the one proposed by Murphy and van der Vaart \cite{MV00} which
uses an asymptotic expansion of approximate profile likelihood.
We summarize the results from the paper.

Suppose we have a function $\eta_{\theta,F}$ that depends on $(\theta
,F)$ such that
$\tilde{\ell}_0(x)\equiv\tilde{\ell}_{\theta_0,F_0}(x)$
is the efficient score function, where
%
\begin{equation}
\label{effscoree} \tilde{\ell}_{\theta,F}(x) \equiv\frac{\partial}{\partial\theta}\log
p_{\theta,\eta_{\theta,F}}(x).
\end{equation}

The theorem below show that if the solution $\hat{\theta}_n$ to the
estimating equation
%
\begin{equation}
\label{EstEqn}
\int\tilde{\ell}_{\hat{\theta}_n,F_n}(x)\,\mathrm{d}F_n=0
\end{equation}
is consistent then it is asymptotically linear with the efficient
influence function
$\tilde{I}_0^{-1}\tilde{\ell}_0(x)$
so that
%
\begin{equation}
\label{AsymptoticNormality}
n^{-1/2}(\hat{\theta}_n -\theta_0)=
\int\tilde{I}_0^{-1}\tilde{\ell }_0(x) \,\mathrm{d}\bigl
\{n^{-1/2}(F_n-F_0)\bigr\} + \mathrm{o}_P(1)
\stackrel{d} {\longrightarrow} N\bigl(0,\tilde{I}_0^{-1}
\bigr),
\end{equation}
%
where $N(0,\tilde{I}_0^{-1})$ is a normal distribution with mean zero
and variance $\tilde{I}_0^{-1}$.
Since $\tilde{I}_0=E_0(\tilde{\ell}_0\tilde{\ell}_0^T)$ is the
efficient information matrix,
this demonstrates that the estimator $\hat{\theta}_n$ is efficient.

On the set of c.d.f. functions ${\cal F}$, we use the sup-norm, that is,
for $F, F_0 \in{\cal F}$,
\[
\|F-F_0\|=\sup_x\bigl|F(x)-F_0(x)\bigr|.
\]
For $\rho>0$, let
\[
{\cal C}_{\rho}=\bigl\{F \in{\cal F}\dvt\|F-F_0\| < \rho\bigr\}.
\]



\begin{thm}[(Hirose \cite{Hirose10})]\label{thmprofile}
Assumptions:
\begin{enumerate}[(R3)]
\item[(R0)]
The function $g_{\theta,F}$ satisfies $g_{\theta_0,F_0}=g_0$ and the function
\[
\tilde{\ell}_0(x)=\tilde{\ell}_{\theta_0,F_0}(x)
\]
is the efficient score function where $\tilde{\ell}_{\theta,F}(x)$ is
given by (\ref{effscoree}).

\item[(R1)] The empirical process $F_n$ is $n^{1/2}$-consistent, that
is, $n^{1/2}\|F_n-F_0\|=\mathrm{O}_P(1)$, and
there exists a $\rho>0$ and a neighborhood $\Theta$ of $\theta_0$
such that for each $(\theta,F) \in\Theta\times{\cal C}_{\rho}$,
the log-likelihood function $\log p(x;\theta,\hat{g}_{\theta,F})$
is twice continuously differentiable with respect to~$\theta$ and
Hadamard differentiable with\vspace*{1pt} respect to $F$ for all $x$.

\item[(R2)] The efficient information matrix $\tilde{I}_0=E_0(\tilde
{\ell}_0\tilde{\ell}_0^{T})$ is invertible.
\item[(R3)] There exists a $\rho>0$ and a neighborhood $\Theta$ of
$\theta_0$
such that
the class of functions $\{\tilde{\ell}_{\theta,F}(x)\dvt (\theta,F) \in
\Theta\times{\cal C}_{\rho}\}$
is Donsker with square integrable envelope function, and that
the class of functions $\{\frac{\partial}{\partial\theta}\tilde{\ell
}_{\theta,F}(x)\dvt (\theta,F) \in\Theta\times{\cal C}_{\rho}\}$
is Glivenko--Cantelli with integrable envelope function.
\end{enumerate}

Under the assumptions $\{(\mathrm{R}0),(\mathrm{R}1),(\mathrm{R}2),(\mathrm{R}3)\}$,
for a consistent solution $\hat{\theta}_n$ to the estimating equation
(\ref{EstEqn}),
the equation (\ref{AsymptoticNormality}) holds.

\end{thm}

\subsection{Asymptotic normality and efficiency in Example~2}\label{sec61}

In this section, we demonstrate how the result of the paper can be used to
show the efficiency of profile likelihood estimators in semi-parametric models.
We show the efficiency of the estimator in Example~2 (using the result
in Section~\ref{sec5}).
First, we identify the efficient score function in the example. Then we
verify conditions (R0)--(R3) in Theorem~\ref{thmprofile}.
Then the efficiency of the estimator follows from the theorem.

\textit{Efficient score function}.
We show that the function (\ref{hatg}) (the solution to the equation
(\ref{g=Psig})) gives us the efficient score function in Example~2.
The log-density function in Example~2 is given by
%
\begin{equation}
\label{logpex2}
\log p(s,z;\theta,g)=1_{\{s=1\}}\bigl\{\log f(y|x;\theta)+
\log g(x)\bigr\} + 1_{\{
s=2\}}\log f_{Y}(y;\theta,g),
\end{equation}
where $z=(y,x)$ if $s=1$ and $z=y$ if $s=2$, and $f_{Y}(y;\theta,g)$ is
given in (\ref{fY}).
%
\begin{thm}[(The efficient score function)]\label{effscoret}
Let us denote $g_{\theta,F_0}(x)$ as the function (\ref{hatg})
evaluated at $(\theta,F_0)$:
%
\begin{equation}
\label{gthetaF0}
g_{\theta,F_0}(x)=\frac{ w_{10} \partial_x \int \mathrm{d}F_{10} }{1 - w_{20}
\int {f(y|x;\theta)}/{f_Y(y;\theta,g_{\theta,F_0})} \,\mathrm{d}F_{20}}.
\end{equation}
%
Then the function
%
\begin{equation}
\label{effscore2}
\tilde{\ell}_{\theta_0,F_0}(s,z) = \frac{\partial}{\partial\theta}
\bigg|_{\theta=\theta_0}\log p(s,z;\theta,g_{\theta,F_0})
\end{equation}
is the efficient score function in the model in Example~2.
\end{thm}

\begin{pf}
We check conditions (\ref{eqnEffCond1}) and (\ref{eqnEffCond2}) in
Theorem~\ref{thmA1} in the \hyperref[appA]{Appendix}.
Then the claim follows from the theorem.

Condition (\ref{eqnEffCond1}) is checked in equation (\ref{VerifyA2ex2}).

We verify condition (\ref{eqnEffCond2}).
Let $g_t(x)$ be a path in the space of density functions with
$g_{t=0}(x)=g_0(x)$.
Define $\alpha_t(x)=g_t(x)-g_0(x)$ and write $\dot{\alpha}_0(x)=\frac
{\partial}{\partial t}|_{t=0} \alpha_t(x)$.
Then
\begin{eqnarray*}
& & \frac{\partial}{\partial t} \bigg|_{t=0}\int\log p(s,z;\theta ,g_{\theta,F_0}+
\alpha_t) \,\mathrm{d}F_0
\\
&&\quad =  \frac{\partial}{\partial t} \bigg|_{t=0} \biggl[ w_{10}\int\bigl\{
\log f(y|x;\theta) + \log(g_{\theta,F_0}+\alpha_t)\bigr
\}\,\mathrm{d}F_{10}\\
&&\qquad \hspace*{31pt}{}+ w_{20} \int \log f_Y(y;
\theta,g_{\theta,F_0}+\alpha_t) \,\mathrm{d}F_{20} \biggr]
\\
&&\quad  =  w_{10}\int\frac{\dot{\alpha}_0(x)}{g_{\theta,F_0}(x)}\,\mathrm{d}F_{10} + w_{20}
\int\frac{\int f(y|x;\theta)\dot{\alpha}_0(x)\,\mathrm{d}x }{f_Y(y;\theta
,g_{\theta,F_0})} \,\mathrm{d}F_{20}
\\
&& \quad =  \int\dot{\alpha}_0(x)\,\mathrm{d}x = \frac{\partial}{\partial t} \bigg|_{t=0}
\int g_t(x) \,\mathrm{d}x=0\quad \bigl(\mbox{by (\protect\ref{gthetaF0}) and since }g_t(x) \mbox{ is a density}\bigr).
\end{eqnarray*}
\upqed\end{pf}

\textit{Efficiency of the profile likelihood estimator}.
Let $\tilde{\ell}_{\theta,F}(s,x)$ be the score function given by (\ref
{effscore2}) with $\theta_0$ and $F_0$ are replaced by $\theta$ and $F$.

We verify conditions $(\mathrm{R}0)$, $(\mathrm{R}1)$, $(\mathrm{R}2)$ and $(\mathrm{R}3)$ of Theorem~\ref{thmprofile} so that we can apply the theorem
to show that the solution $\hat{\theta}_n$ to the estimating equation
\[
\sum_{s=1}^2\sum
_{i=1}^n\tilde{\ell}_{\hat{\theta}_n,F_n}(s,X_{si})=0
\]
is asymptotically linear estimator with the efficient influence
function, that is, (\ref{AsymptoticNormality}) holds.
This shows the efficiency of the MLE based on the profile likelihood in
this example.

\textit{Condition (R0)}.
Theorem~\ref{effscoret}  shows
that the score function evaluated at $(\theta_0,F_0)$ is the efficient
score function in Example~2.

\textit{Condition (R1)}.
We assume that:
\begin{enumerate}[(T1)]
\item[(T1)] For all $\theta\in\Theta$, the function $f(y|x;\theta)$
is twice continuously differentiable with respect to~$\theta$.
\end{enumerate}

The maps
\[
g \rightarrow\log g(x)
\]
and
\[
g \rightarrow f_Y(y;\theta,g)=\int_{\cal X} f(y|x;
\theta)g(x)\,\mathrm{d}x
\]
are Hadamard differentiable (cf. Gill \cite{Gill89}).
It follows that the log-density function
$\log p(s,z;\theta,g)$ given by (\ref{logpex2})
is Hadamard differentiable with respect to $g$ and, by assumption (T1),
it is also twice continuously differentiable with respect to $\theta$.
In the previous section (Section~\ref{sec5}), we verified the function
$g_{\theta,F}$ is Hadamard differentiable with respect to $F$ and twice
continuously differentiable with respect to $\theta$.
By the chain rule and product rule of Hadamard differentiable maps, the
log-density function
$\log p(s,x;\theta,g_{\theta,F})$ is
Hadamard differentiable with respect to $F$ and twice continuously
differentiable with respect to $\theta$.
Therefore, we verified condition (R1).

\textit{Derivatives of log-likelihood}.
The log-density function under consideration is
%
\begin{equation}
\label{EqnLogLikelihoodOneInduced}
\log p(s,z;\theta,g_{\theta,F})  =  1_{\{s=1\}}\bigl\{
\log f(y|x;\theta )+\log g_{\theta,F}(x)\bigr\} + 1_{\{s=2\}}\log
f_{Y}(y;\theta,g_{\theta
,F}).
\end{equation}
The derivative of the log-density with respect to $\theta$ is
%
\begin{eqnarray}
\tilde{\ell}_{\theta,F}(s,z) & = & \frac{\partial}{\partial\theta} \log
p(s,z;\theta,g_{\theta,F})
\nonumber
\\[-8pt]
\label{EqnScorefunction}\\[-8pt]
\nonumber
& = & 1_{\{s =1\}} \biggl\{ \frac{\dot{f}}{f} + \frac{\dot{g}_{\theta
,F}}{g_{\theta,F}} \biggr
\} + 1_{\{s =2\}} \frac{\dot{f}_Y+ \mathrm{d}_g f_Y(\dot{g}_{\theta,F})}{f_Y}.
\end{eqnarray}

The second derivative of the log-density function with respect to
$\theta$ is
%
\begin{eqnarray}
\frac{\partial}{\partial\theta^T}\tilde{\ell}_{\theta,F}(s,z) &=&
\frac{\partial^2}{\partial\theta \,\partial\theta^T} \log p(s,z;\theta ,g_{\theta,F})
\nonumber
\\
& = & 1_{\{s =1\}} \biggl\{ \frac{\ddot{f}}{f}-\frac{\dot{f}\dot{f}^T}{f^2} +
\frac{\ddot{g}_{\theta,F}}{g_{\theta,F}}-\frac{\dot{g}_{\theta,F}\dot
{g}^T_{\theta,F}}{g^2_{\theta,F}} \biggr\}
\nonumber
\\
\label{EqnScoreDerivTheta}
&& {}+ 1_{\{s =2\}} \biggl\{ \frac{\ddot{f}_Y+ \mathrm{d}_g \dot{f}_Y(\dot
{g}_{\theta,F})}{f_Y} -\frac{\dot{f}_Y\dot{f}^T_Y+ \dot{f}_Y \,\mathrm{d}_g f_Y(\dot{g}^T_{\theta
,F})}{f^2_Y}
\\
&& \hspace*{40pt}{}+ \frac{ \mathrm{d}_g \dot{f}^T_Y(\dot{g}_{\theta
,F})+\mathrm{d}_g f_Y(\ddot{g}_{\theta,F})}{f_Y} \nonumber\\
&&\hspace*{40pt}{}-\frac{\mathrm{d}_g f_Y(\dot{g}_{\theta,F})\dot{f}^T_Y+ \mathrm{d}_g f_Y(\dot{g}_{\theta
,F}) \,\mathrm{d}_g f_Y(\dot{g}^T_{\theta,F})}{f^2_Y} \biggr\}.\nonumber
\end{eqnarray}
Here, we used the notation
$\dot{f}_Y=\dot{f}_Y(y;\theta,g_{\theta,F})$, $\ddot{f}_Y=\ddot
{f}_Y(y;\theta,g_{\theta,F})$,
$\mathrm{d}_g f_Y(g_{\theta,F})=\int f(y|x;\theta)g_{\theta,F}(x)\,\mathrm{d}x$, and $\mathrm{d}_g
\dot{f}_Y(g_{\theta,F})=\int\dot{f}(y|x;\theta)g_{\theta,F}(x)\,\mathrm{d}x$.

\textit{Condition (R2)}.
We assume that:
\begin{enumerate}
\item[(T2)] There is no $a \in R^d$ such that $a^T\frac{\dot
{f}}{f}(y|x;\theta)$ is
constant in $y$ for almost all $x$.
\end{enumerate}
The term $\frac{\dot{g}_{\theta,F}}{g_{\theta,F}}(x,\theta_0,F_0)$ is a
function of $x$.
Therefore, by equation~(\ref{EqnScorefunction}) and assumption (T2),
there is no $a \in R^d$ such that $a^T \tilde{\ell}_{\theta,F}(1,z)$ is
constant in $y$ for almost all $x$.
By Theorem~1.4 in Seber and Lee \cite{SL03},
$E (\tilde{\ell}_{\theta_0,F_0} \tilde{\ell}^T_{\theta_0,F_0})$ is
nonsingular with the bounded inverse.

\textit{Conditions (R3)}.
Since verification of condition (R3) require more assumptions and it
does not add anything new, we simply assume:
\begin{enumerate}[(T3)]
\item[(T3)]
Let ${\cal F}$ be the set of c.d.f. functions and for some $\rho>0$ define
${\cal C}_{\rho}=\{F \in{\cal F}\dvt \|F-F_0\|_{\infty}\leq\rho\}$.
The class of function
\[
\bigl\{\tilde{\ell}_{\theta,F}(s,z)\dvt (\theta,F) \in\Theta \times {\cal
C}_{\rho} \bigr\}
\]
is $P_{\theta_0,g_0}$-Donsker with square integrable envelope function
and the class
\[
\biggl\{\frac{\partial}{\partial\theta^T}\tilde{\ell}_{\theta,F}(s,z)\dvt (\theta,F) \in
\Theta \times{\cal C}_{\rho} \biggr\}
\]
is $P_{\theta_0,g_0}$-Glivenko--Cantelli with integrable envelope function.
\end{enumerate}

\section{Discussion}\label{sec7}
In Theorem~\ref{mainthm}, we have shown the differentiability
of implicitly defined function which we encounter in the maximum
likelihood estimation in semi-parametric models.
In the theorem, we assumed the implicitly defined function is the
solution to the operator equation (\ref{ImplicitEquation})
and we obtained the derivatives of the (implicitly defined) function.
In application of the theorem,~we need to verify condition (A3) in the
theorem (that is $\|\mathrm{d}_{\eta}\Psi_{\theta_0,F_0}(\eta_0)\|<1$).
This required additional conditions in the examples ((\ref
{ConditionA3}) in Example~1 and (\ref{ConditionA3ex2}) in
Example~2).
The future work is to relax the condition to $\|\mathrm{d}_{\eta}\Psi_{\theta
_0,F_0}(\eta_0)\|< \infty$ so that the additional conditions can be weaken.
Once the differentiability of the implicitly defined function has been
established,
the results in Hirose \cite{Hirose10} (we summarized in Section~\ref{sec6},
Theorem~\ref{thmprofile}) are applicable.


\begin{appendix}

\section{Verification of efficient score function}\label{appA}

To verify condition $(\mathrm{R}0)$ in Theorem~\ref{thmprofile}, the
following theorem may be useful.
This is a modification of the proof in Breslow, McNeney and Wellner \cite{BMW00} which was
originally adapted from Newey \cite{Newey94}.

\begin{thm}\label{thmA1}
We assume the general semi-parametric model given in  the \hyperref[sec1]{Introduction}
with the density $p_{\theta,\eta}(x)=p(x;\theta,\eta)$ is
differentiable with respect to $\theta$
and Hadamard differentiable with respect to $\eta$.
Suppose $g_t$ is an arbitrary path such that $g_{t=0}=g_0$ and let
$\alpha_t=g_t-g_0$.
If $g_{\theta,F}$ is a function of $(\theta,F)$ such that
%
\begin{equation}
\label{eqnEffCond1}
g_{\theta_0,F_0}=g_0
\end{equation}
and, for each $\theta\in\Theta$,
%
\begin{equation}
\label{eqnEffCond2}
\frac{\partial}{\partial t} \bigg|_{t=0} E_0 \bigl[\log
p(x;\theta ,g_{\theta,F_0}+\alpha_t) \bigr]=0,
\end{equation}
then the function $\tilde{\ell}_{\theta_0,F_0}(x)=\frac{\partial
}{\partial\theta} |_{\theta=\theta_0} \log p(x;\theta,g_{\theta,F_0})$
is the efficient score function.
\end{thm}

\begin{pf}
Condition (\ref{eqnEffCond2}) implies that
%
\begin{eqnarray}
0 & = & \frac{\partial}{\partial\theta} \bigg|_{\theta=\theta_0} \frac
{\partial}{\partial t}
\bigg|_{t=0} E_{0} \bigl[\log p(x;\theta,g_{\theta
,F_0}+
\alpha_t) \bigr]
\nonumber
\\[-8pt]
\label{eqnEffCond3}\\[-8pt]
\nonumber
& = & \frac{\partial}{\partial t} \bigg|_{t=0} E_{0} \biggl[
\frac{\partial
}{\partial\theta} \bigg|_{\theta=\theta_0}\log p(x;\theta,g_{\theta
,F_0}+
\alpha_t) \biggr].
\end{eqnarray}

By differentiating the identity
\[
\int \biggl(\frac{\partial}{\partial\theta}\log p(x;\theta,g_{\beta
,F_0}+
\alpha_t) \biggr)p(x;\theta,g_{\beta,F_0}+\alpha_t)\,\mathrm{d}x=0
\]
with respect to $t$ at $t=0$ and $\theta=\theta_0$, we get
%
\begin{eqnarray}
0 & = & \frac{\partial}{\partial t} \bigg|_{t=0,\theta=\theta_0}\int \biggl(
\frac{\partial}{\partial\theta}\log p(x;\theta,g_{\theta,F_0}+\alpha _t)
\biggr)p(x;\theta,g_{\theta,F_0}+\alpha_t)\,\mathrm{d}x
\nonumber
\\
& = & E_{0} \biggl[\tilde{\ell}_{\theta_0,F_0}(x) \biggl(
\frac{\partial
}{\partial t} \bigg|_{t=0}\log p(x;\theta_0,g_t)
\biggr) \biggr] \qquad  \bigl(\mbox{by (\ref{eqnEffCond1})}\bigr)
\nonumber
\\[-8pt]
\label{eqnEffCond4}\\[-8pt]
\nonumber
&&{}+\frac{\partial}{\partial t} \bigg|_{t=0} E_{0} \biggl[
\frac
{\partial}{\partial\theta} \bigg|_{\theta=\theta_0}\log p(x;\theta ,g_{\theta,F_0}+
\alpha_t) \biggr]
\nonumber
\\
& = & E_{0} \biggl[\tilde{\ell}_{\theta_0,F_0}(x) \biggl(
\frac{\partial
}{\partial t} \bigg|_{t=0}\log p(x;\theta_0,g_t)
\biggr) \biggr] \qquad  \bigl(\mbox{by (\ref{eqnEffCond3})}\bigr).\nonumber
\end{eqnarray}

Let $c \in R^m$ be arbitrary.
Then it follows from equation (\ref{eqnEffCond4}) that the product
$c'\tilde{\ell}_{\theta_0,F_0}(x)$ is orthogonal to the nuisance
tangent space $\dot{\cal P}_g$
which is the closed linear span of score functions of the form $\frac
{\partial}{\partial t} |_{t=0}\log p(x;\beta_0,g_t)$.

Using condition (\ref{eqnEffCond1}), we have
\begin{eqnarray*}
\tilde{\ell}_{\theta_0,F_0}(x) & = & \frac{\partial}{\partial\theta} \bigg|_{\theta=\theta_0} \log
p(x;\theta,g_0)+ \frac{\partial}{\partial\beta} \bigg|_{\theta=\theta_0} \log p(x;\theta
_0,g_{\theta,F_0})
\\
& = & \dot{\ell}_{\theta_0,g_0}(x)- \psi_{\theta_0,g_0}(x),
\end{eqnarray*}
where $\dot{\ell}_{\theta_0,g_0}(x)=\frac{\partial}{\partial\theta}
|_{\theta=\theta_0} \log p(x;\theta,g_0)$ is the score function for
$\theta$
and
$\psi_{\theta_0,g_0}(x)= -\frac{\partial}{\partial\theta} |_{\theta=\theta_0} \log p(x;\theta_0,g_{\theta,F_0})$.
Finally, $c'\tilde{\ell}_{\theta_0,F_0}(x)=c'\dot{\ell}_{\theta
_0,g_0}(x)-c'\psi_{\theta_0,g_0}(x)$ is orthogonal to the nuisance
tangent space $\dot{\cal P}_{g}$ and $c'\psi_{\theta_0,g_0}(x) \in\dot
{\cal P}_{g}$ implies that
$c'\psi_{\theta_0,g_0}(x)$ is the orthogonal projection of $c'\dot{\ell
}_{\theta_0,g_0}(x)$ onto the nuisance tangent space $\dot{\cal P}_{g}$.
Since $c \in R^m$ is arbitrary, $\tilde{\ell}_{\theta_0,F_0}(x)$ is the
efficient score function.
\end{pf}

\section{Verification of (A1) in Example~1: Continued from Section~\texorpdfstring{\protect\ref{sec4}}{4}}\label{appB}

In verification of (A1) in Example~1, Section~\ref{sec4}, we gave proof
the Hadamard differentiability of functions with additional condition
for the derivative of $\Psi_{\beta,F}(A)$ with respect to $F$.
For the rest the derivatives, we give them without proofs.

(\textit{The derivative of $\Psi_{\beta,F}(A)$ with respect to $A$})
Let $h_1=h_1(U)$ be a function of $U$.
%
\begin{equation}
\label{dAPsi}
\mathrm{d}_A\Psi_{\beta,F}(A)h_1 =
-E_{F}\int_0^u
\frac{E_{F} \,\mathrm{d}_AW(s;\beta
,A)h_1 }{\{E_{F} W(s;\beta,A)\}^2 }\,\mathrm{d}N(s),
\end{equation}
where
%
\begin{equation}
\label{dAW} \mathrm{d}_A W(s;\beta,A)h_1= \frac{-(1+\delta)\mathrm{e}^{2\beta'Z}Y(s) h_1(U)}{\{
1+\mathrm{e}^{\beta'Z}A(U)\}^2}.
\end{equation}

(\textit{The second derivative of $\Psi_{\beta,F}(A)$ with respect to $A$})
If $h_1(U)$, $h_2(U)$ are functions,
\begin{eqnarray*}
\mathrm{d}^2_A\Psi_{\beta,F}(A)h_1h_2
& = & E_{F}\int_0^u
\frac{E_{F}\,
\mathrm{d}^2_AW(s;\beta,A)h_1h_2 }{\{E_{F} W(s;\beta,A)\}^2 }\,\mathrm{d}N(s)
\\
&&{}+ E_{F}\int_0^u
\frac{2 \{E_{F}\, \mathrm{d}_AW(s;\beta,A)h_1\}\{ E_{F}\,\mathrm{d}_A
W(s;\beta,A)h_2\}
}{ \{E_{F} W(s;\beta,A)\}^3}\, \mathrm{d}N(s),
\end{eqnarray*}
where
\[
\mathrm{d}^2_A W(s;\beta,A_t)h_1h_2=
\frac{2(1+\delta)\mathrm{e}^{3\beta'Z}Y(s)
h_1(U)h_2(U) }{\{1+\mathrm{e}^{\beta'Z}A(U)\}^3}.
\]
(The expression of $\mathrm{d}_AW(s;\beta,A)h_1$ is given in (\ref{dAW}).)

(\textit{The first and second derivative of $\Psi_{\beta,F}(A)$ with
respect to $\beta$})
Let us denote the first and second derivatives (with respect to $\beta
$) by $\dot{\Psi}_{\beta,F}(A)$
and $\ddot{\Psi}_{\beta,F}(A)$, respectively. Then they are given by,
for $a, b \in R^d$,
\begin{eqnarray*}
a^T \dot{\Psi}_{\beta,F}(A) & = & a^T \biggl\{
\frac{\partial}{\partial
\beta}\Psi_{\beta,F}(A) \biggr\} \\
&= & -E_{F}\int
_0^u \frac{E_{F} a^T \dot
{W}(s;\beta,A) }{\{E_{F} W(s;\beta,A)\}^2 }\,\mathrm{d}N(s),
\\
a^T \ddot{\Psi}_{\beta,A}(g)b & = & a^T \biggl\{
\frac{\partial^2}{\partial\beta\, \partial\beta^T}\Psi _{\beta,F}(A) \biggr\}b
\\
& = & E_{F}\int_0^u
\frac{E_{F} a^T\ddot{W}(s;\beta,A)b }{\{E_{F}
W(s;\beta,A)\}^2 }\,\mathrm{d}N(s)
\\
& & + E_{F}\int_0^u
\frac{2 \{E_{F} a^T\dot{W}(s;\beta,A)h_1\}\{
E_{F}\dot{W}^T(s;\beta,A)b\}
}{ \{E_{F} W(s;\beta,A)\}^3} \,\mathrm{d}N(s).
\end{eqnarray*}
Here,
\[
a^T\dot{W}(s;\beta,A)  =  a^T \biggl\{
\frac{\partial}{\partial\beta
}W(s;\beta,A) \biggr\} 
=
\frac{(1+\delta)a^T\beta \mathrm{e}^{\beta^TZ}Y(s)}{\{1+\mathrm{e}^{\beta^TZ}A(U)\}^2}
\]
and
\begin{eqnarray*}
a^T\ddot{W}(s;\beta,A)b  &=&  a^T \biggl\{
\frac{\partial^2}{\partial
\beta\,\partial\beta^T}W(s;\beta,A) \biggr\}b
\\
& = &
\frac{(1+\delta)\{(a^Tb)\mathrm{e}^{\beta^TZ}+(a^T\beta)(\beta^Tb) \mathrm{e}^{\beta
^TZ}\}Y(s)}{\{1+\mathrm{e}^{\beta^TZ}A(U)\}^2} \\
&&{}-\frac{2(1+\delta)(a^T\beta)(\beta^Tb) \mathrm{e}^{2\beta^TZ}Y(s)A(U)}{\{
1+\mathrm{e}^{\beta^TZ}A(U)\}^3}.
\end{eqnarray*}

(\textit{The derivative of $\Psi_{\beta,F}(A)$ with respect to $\beta$ and $A$})
For given function $h_1(U)$ and $a \in R^d$,
\begin{eqnarray*}
&& a^T\, \mathrm{d}_A \dot{\Psi}_{\beta,F}(A)h_1
\\
&&\quad  =- E_{F}\int_0^u \biggl\{
\frac{E_{F} a^T \,\mathrm{d}_A\dot{W}(s;\beta,A)h_1
}{\{E_{F} W(s;\beta,A)\}^2 } -2 \frac{E_{F} a^T \dot{W}(s;\beta,A)
E_{F}\, \mathrm{d}_A W(s;\beta,A)h_1 }{\{E_{F} W(s;\beta,A)\}^3 } \biggr\}\,\mathrm{d}N(s),
\end{eqnarray*}
here $ a^T \dot{W}(s;\beta,A)$ is given above, $\mathrm{d}_A W(s;\beta,A)h_1$
is given in (\ref{dAW}) and
\[
a^T\,\mathrm{d}_A \dot{W}(s;\beta,A)h_1=
\frac{-2(1+\delta)a^T\beta \mathrm{e}^{2\beta
^TZ}Y(s) h_1(U)}{\{1+\mathrm{e}^{\beta^TZ}A(U)\}^3}.
\]

\section{Verification of (A1) in Example~2: Continued from~Section~\texorpdfstring{\protect\ref{sec5}}{5}}\label{appC}

We proved the Hadamard differentiability of functions and additional
condition for the derivative of $\Psi_{\theta,F}(g)$ with respect to
$F$ in Section~\ref{sec5},
verification of (A1) in Example~2.
The rest of the derivatives are listed here.

(\textit{The derivative of $\Psi_{\theta,F}(g)$ with respect to $g$})
For a function $h^*(x)$ of $x$,
%
\begin{equation}
\label{dgPsi} \mathrm{d}_g\Psi_{\theta,F}(g)h^* =  \frac{- (\partial_x \int\pi_1(\mathrm{d}F) ) \{\mathrm{d}_gA(x;\theta
,g,F)h^* \}}{\{A(x;\theta,g,F)\}^2},
\end{equation}
where
%
\begin{equation}
\label{dgA}
\mathrm{d}_gA(x;\theta,g,F)h^*  =  \int f(y|x;\theta)
\frac{\int f(y|x;\theta
)h^*(x)\,\mathrm{d}x }{\{f_Y(y;\theta,g)\}^2} \pi_2(\mathrm{d}F).
\end{equation}

(\textit{The second derivative of $\Psi_{\theta,F}(g)$ with respect to $g$})
For functions $h_1(x)$ and $h_2(x)$ of $x$,
\begin{eqnarray*}
\label{d^2gPsi}
&&
\mathrm{d}^2_g\Psi_{\theta,F}(g)h_1h_2
\\
&&\quad  =  \biggl(\partial_x \int\pi_1(\mathrm{d}F) \biggr) \biggl[-
\frac{\mathrm{d}^2_gA(x;\theta
,g,F)h_1h_2}{\{A(x;\theta,g_t,F)\}^2} +\frac{2\{\mathrm{d}_gA(x;\theta,g,F)h_1\}\{\mathrm{d}_gA(x;\theta,g,F)h_2\}}{\{
A(x;\theta,g,F)\}^3} \biggr],
\end{eqnarray*}
where
\[
\mathrm{d}^2_gA(x;\theta,g_t,F)h_1h_2
 =  -2\int f(y|x;\theta) \frac{ \{\int f(y|x;\theta
)h_1(x)\,\mathrm{d}x \}  \{\int f(y|x;\theta)h_2(x)\,\mathrm{d}x \} }{\{
f_Y(y;\theta,g)\}^3} \pi_2(\mathrm{d}F).
\]

(\textit{The first and second derivative of $\Psi_{\theta,F}(g)$ with
respect to $\theta$})
Let us denote the first and second derivatives with respect to $\theta$
by $\dot{\Psi}_{\theta,F}(g)$
and $\ddot{\Psi}_{\theta,F}(g)$, respectively. They are given by, for
$a, b \in R^d$,
\begin{eqnarray*}
\label{dotPsi} a^T\dot{\Psi}_{\theta,F}(g) & = & a^T
\biggl\{ \frac{\partial}{\partial
\theta}\Psi_{\theta,F}(g) \biggr\} = -
\frac{ (\partial_x \int\pi_1(\mathrm{d}F) ) a^T\dot{A}}{A^2},
\\
\label{ddotPsi} a^T \ddot{\Psi}_{\theta,F}(g)b & = &
a^T \biggl\{ \frac{\partial
^2}{\partial\theta\, \partial\theta^T}\Psi_{\theta,F}(g) \biggr\}b = -
\frac{ (\partial_x \int\pi_1(\mathrm{d}F) )\{A(a^T \ddot
{A}b)-2(a^T\dot{A})(\dot{A}^T b)\}}{A^3},
\end{eqnarray*}
where
\[
a^T\dot{A}  =   a^T \biggl\{\frac{\partial}{\partial\theta}A(x;
\theta ,g,F) \biggr\} = - \int \frac{ f_Y (a^T\dot{f})-f(a^T\dot{f}_Y)}{f_Y^2} \pi_2(\mathrm{d}F)
\]
and
\begin{eqnarray*}
a^T\ddot{A}b  &=&  a^T \biggl\{\frac{\partial^2}{\partial\theta\,
\partial\theta^T}A(x;
\theta,g,F) \biggr\}b
\\
& = & - \int \bigl( f_Y^2\bigl(a^T \ddot{f} b\bigr)-ff_Y\bigl(a^T \ddot{f}_Y
b\bigr)+2f\bigl(a^T\dot{f}_Y\bigr)\bigl(\dot{f}_Y^Tb\bigr)\\
&&\hspace*{23pt}{}-f_Y\bigl(a^T\dot{f}\bigr)\bigl(\dot{f}_Y^Tb\bigr)
-f_Y\bigl(a^T\dot{f}_Y\bigr)\bigl(\dot{f}^Tb\bigr)\bigr)/{f_Y^3} \pi_2(\mathrm{d}F).
\end{eqnarray*}

(\textit{The derivative of $\Psi_{\theta,F}(g)$ with respect to $\theta$
and $g$})
For $a \in R^d$ and function $h^*(x)$ of $x$,
\begin{eqnarray*}
\label{dgdotPsi}
&& a^T \,\mathrm{d}_g \dot{\Psi}_{\theta,F}(g)h^*
\\
&& \quad =  - \biggl( \partial_x \int\pi_1(\mathrm{d}F) \biggr) \biggl[
\frac{a^T \,\mathrm{d}_g\dot{A}(x;\theta,g,F)h^*}{\{A(x;\theta,g,F)\}^2} - \frac{2 a^T \dot{A}(x;\theta,g,F) \,\mathrm{d}_g A(x;\theta,g,F)h^*}{\{A(x;\theta
,g,F)\}^3} \biggr],
\end{eqnarray*}
where
\begin{eqnarray*}
\label{dgdotA} &&\hspace*{-3pt} a^T \,\mathrm{d}_g\dot{A}(x;\theta,g,F)h^*
\\
&&
\hspace*{-5pt}\quad =  \int\bigl(a^T \dot{f}\bigr) \frac{\int f h^*\,\mathrm{d}x}{f_Y^2}
\pi_2(\mathrm{d}F) + \int f \frac{\int(a^T\dot{f})h^*\,\mathrm{d}x}{f_Y^2} \pi_2(\mathrm{d}F) - 2\int f
\bigl(a^T\dot{f}_Y\bigr)\frac{ \int fh^*\,\mathrm{d}x }{f_Y^3}
\pi_2(\mathrm{d}F).
\end{eqnarray*}
\end{appendix}






\printhistory
\end{document}